\documentclass[1p]{elsarticle}

\usepackage{lineno,hyperref}
\usepackage{graphicx,amssymb,amsmath,amsfonts,amsthm}
\usepackage[utf8]{inputenc}
\usepackage{mathtools}
\usepackage{color}
\usepackage{bbm}
\usepackage{arydshln} 
\usepackage{transparent}
\usepackage{subfigure}
\usepackage{float}
\usepackage{subfloat}
\usepackage{accents}
\hypersetup{backref,pdfpagemode=FullScreen,colorlinks=true,citecolor=red,linkcolor=blue}

\newtheorem{lemma}{Lemma}
\newtheorem{proposition}{Proposition}
\newtheorem{remark}{Remark}

\newtheorem{theorem}{Theorem}
\newtheorem{corollary}{Corollary}
\theoremstyle{definition}
\newtheorem{definition}{Definition}
\theoremstyle{definition}

\modulolinenumbers[5]

\journal{Journal of Geometry and Physics}

\begin{document}

\begin{frontmatter}

\title{Principal cycles of one dimensional foliations associated to a plane field  in   $\mathbbm{E}^3$.}

\author[1]{Alacyr J. Gomes\corref{a}}
\author[2]{Ronaldo A. Garcia\corref{a}}
\address{Universidade Federal de  Goi\'as,
	Instituto de Matem\'atica e Estat\'\i stica.\corref{a}}
\fntext[1]{E-mail address: alacyr@ufg.br}
\fntext[2]{E-mail address: ragarcia@ufg.br}
\cortext[a]{Universidade Federal de  Goi\'as,
	Instituto de Matem\'atica e Estat\'\i stica,
	Campus Samambaia, Esquina da Rua Jacarand\'a com a Av. Pau-Brasil,
	74.690-900,  Goi\^ania, Goi\'as, Brasil.}

\begin{abstract}
	In this work  it  will be  analyzed $\eta$-principal cycles (compact leaves) of  one dimensional singular foliations  associated to    a plane field $\Delta_{\eta}$ defined by a unit and normal vector field ${\eta}$ in $ \mathbbm E^3$. The leaves are orthogonal to the orbits of  ${\eta}$  and are the integral curves corresponding to directions of extreme   normal curvature of the plane field  $\Delta_{\eta}$.
It  is shown that, generically, given a $\eta$-principal cycle it  can be  make hyperbolic (the derivative of the first return  of the Poincar\'e map has all eigenvalues disjoint from the unit circle) by a small deformation of the vector field ${\eta}$. Also is shown that for a dense set of unit vector fields, with the weak  $C^r$-topology of Whitney, the $\eta$-principal cycles are hyperbolic.
\end{abstract}

\begin{keyword}
\texttt{vector field \sep principal foliations \sep hyperbolic principal cycle \sep plane field \sep normal curvature.}
\MSC[2010] 37C27 \sep  34C25 \sep 	53C12 \sep 	93B27 \sep 	93C15
\end{keyword}

\end{frontmatter}


\section{Introduction and Main Results}\label{int}

In this paper it will be analyzed the $\eta-$principal configuration defined by two one dimensional foliations associated to a plane field in $\mathbb E^3$  near a compact leaf.  The $\eta-$principal configuration can be described as the extrinsic geometry of  plane fields   in $\mathbb E^3$. It consists of two orthogonal  unidimensional singular foliations which  are orthogonal to the orbits of the vector field defining the plane distribution.

This work was motivated by the extrinsic geometry of vector fields and plane fields  $ \Delta\eta$  as developed by Y. Aminov  \cite{YA-2000}.  This approach came back to the classical works, G. Rogers \cite{ro-1912}, A. Voss \cite{Voss} and others.  See also  \cite{bo-1995},  \cite{Kr-2008} and \cite{Re-1977}.  

Let  $\eta$ be a  unit  smooth vector field  in $\mathbbm{E}^3$ and $\Delta_{\eta}$ be a plane field distribution orthogonal to $\eta.$
The second fundamental form of $\Delta_{\eta}$   is the bilinear form 
\[ II(X,Y)=\frac 12  \langle \nabla_X Y + \nabla_Y X,\eta\rangle\]
where $X$ and $Y$ are smooth vector fields generating $\Delta_{\eta}$, see \cite{Kr-2008} and \cite{Re-1977}.

The normal curvature of $\Delta_{\eta}$ is given by
 \[ k_{\eta}(X)=\frac{II(X,X)}{I(X,X)}, \;\;\;I(X,Y)=\langle X, Y\rangle\]
 
 Following the classical approach of principal curvature lines (see \cite{YA-2000}, \cite{SB-1995},  \cite{CGJS-1982}, \cite{MS03}) the extremal values of $k_{\eta}$ restricted to   the plane field $\Delta_{\eta}$ are called $\eta-$prin\-cipal curvatures (denoted by $k_1\leq k_2$) and the associated directions are called $\eta-$principal directions (denoted by $e_1$ and $e_2$). The integral curves of $e_1$ and $e_2$ are called $\eta-$principal lines, defining the two $\eta-$principal foliations $\mathcal{F}_i(\eta)$, $i=1,2$.
 The closed integral curves are called $\eta-$principal cycles.

The main results are:

\begin{theorem}{\label{th:pert-cont}}
Let $\gamma$  be a $\eta-$principal cycle   of length   $L$
of the $\eta$-principal foliation $\mathcal{F}_1(\eta)$.  Given $\varepsilon>0$ small, there is a smooth vector field $\eta_{\varepsilon}$ in $\mathfrak{X}_{\mathfrak{R}}^r(\mathbbm{E}^3)$,  such that 
$\eta_{\varepsilon} $ is $\varepsilon-C^r $ close to $\eta$,  with $\gamma$ being a  hyperbolic  $\eta_{\varepsilon}-$principal cycle of  $\mathcal{F}_1(\eta_{\varepsilon})$.
\end{theorem}

\begin{theorem}{\label{th:peridoc-dense}}
 Consider the set $\mathcal{G}$ of $\mathfrak{X}_{\mathfrak{R}}^r(\mathbbm{E}^3)$   such that, for $\eta \in \mathcal{G}$, all  $\eta-$principal cycles of  $\mathcal{F}_1(\eta)$ and $\mathcal{F}_2(\eta)$ are hyperbolic.
 Then $\mathcal{G}$ is dense in $\mathfrak{X}_{\mathfrak{R}}^r(\mathbbm{E}^3)$.
\end{theorem}

In section \ref{curv-normal} the basic definitions will be introduced. The main concepts   are  the $\eta-$principal line fields, whose integral curves are   called  $\eta-$principal curvature lines.  They are associated to the normal curvature of  a plane distribution defined by   an unit vector field $\eta$. See  \cite{YA-2000}.   After that, it is obtained the differential equation of the $\eta-$principal line field directions.

When the plane distribution $\Delta_\eta$ is integrable these concepts coincide  with the classical theory of principal lines on surfaces of $\mathbbm E^3$, a classical subject of differential geometry of surfaces which were introduced by G. Monge (1796), \cite{Mo-1796}. The qualitative theory and global aspects of principal lines were initiated  by   C. Gutierrez and J. Sotomayor (1982). See \cite{RGJS-2009}, \cite{JSCG-1991} and \cite{CGJS-1982}.


In Section \ref{LCfechada-ret},  the first return map associated to a $\eta-$principal cycle is considered and its first  derivative is obtained as a solution of a  linear differential equation. It will be   shown that generically  $\eta-$principal cycles are hyperbolic. In order to obtain the main result of Theorem \ref{th:pert-cont}  we will make use of ideas and results of geometric control theory as developed in \cite{RR-2012}  in the context of Hamiltonian systems and geodesics. 

In Sections \ref{sec:th1} and  \ref{sec:th2} will be presented  the proofs of Theorems \ref{th:pert-cont} and  \ref{th:peridoc-dense}.

In Section \ref{sec:exe} examples of hyperbolic  $\eta$-principal cycles are analyzed.

In  Section  \ref{sec:cr} some concluding remarks are discussed.

In \ref{app:A}  the main theorem involving geometric control theory, that will be used in Section \ref{sec:th1}  is stated for completeness and convenience to the reader.

In \ref{app:B} the topology of $\mathfrak{X}_{\mathfrak{R}}^r(\mathbbm{E}^3)$ will be reviewed.


\section{ $\eta-$Principal lines associated to a plane field $\Delta_\eta$ }\label{curv-normal}

In this section,  we present some   notions concerning the extrinsic geometry of vector fields and plane fields and results that will be used in this work. This section is inspired in the Aminov's  book   \cite[Chapter 1]{YA-2000}.

\subsection{Normal curvature of a plane field}


Let $\mathfrak{X}_{\mathfrak{R}}^{r}(\mathbbm{E}^3) $ be the set of unit 
regular vector fields  of class $C^r$ in   $\mathbbm{E}^3$ with the $C^r-$topology of Whitney, see    \cite{HM-1976} and \cite{MP-1967}.  

Associated to a vector field $\eta=(\eta_1,\eta_2,\eta_3)\in \mathfrak{X}_{\mathfrak{R}}^{r}(\mathbbm{E}^3)$  we  have a  plane distribution, or a plane field,  which will be denoted by  $\Delta_{\eta}$. For each $p$, $\Delta(p)$ is the plane having $\eta(p)/|\eta(p)|$ as unit normal vector.
The plane distribution $\Delta_{\eta}$ is completely integrable if, and only if, $\langle \rm{curl}{(\eta )} , \eta \rangle=0$, where $\langle \cdot , \cdot \rangle$ denotes the usual inner product in $\mathbbm{E}^3$, see Frobenius  Theorem  \cite[page 192]{MS01}.

\begin{figure}[H]
	\begin{center}
		\includegraphics[scale=0.25]{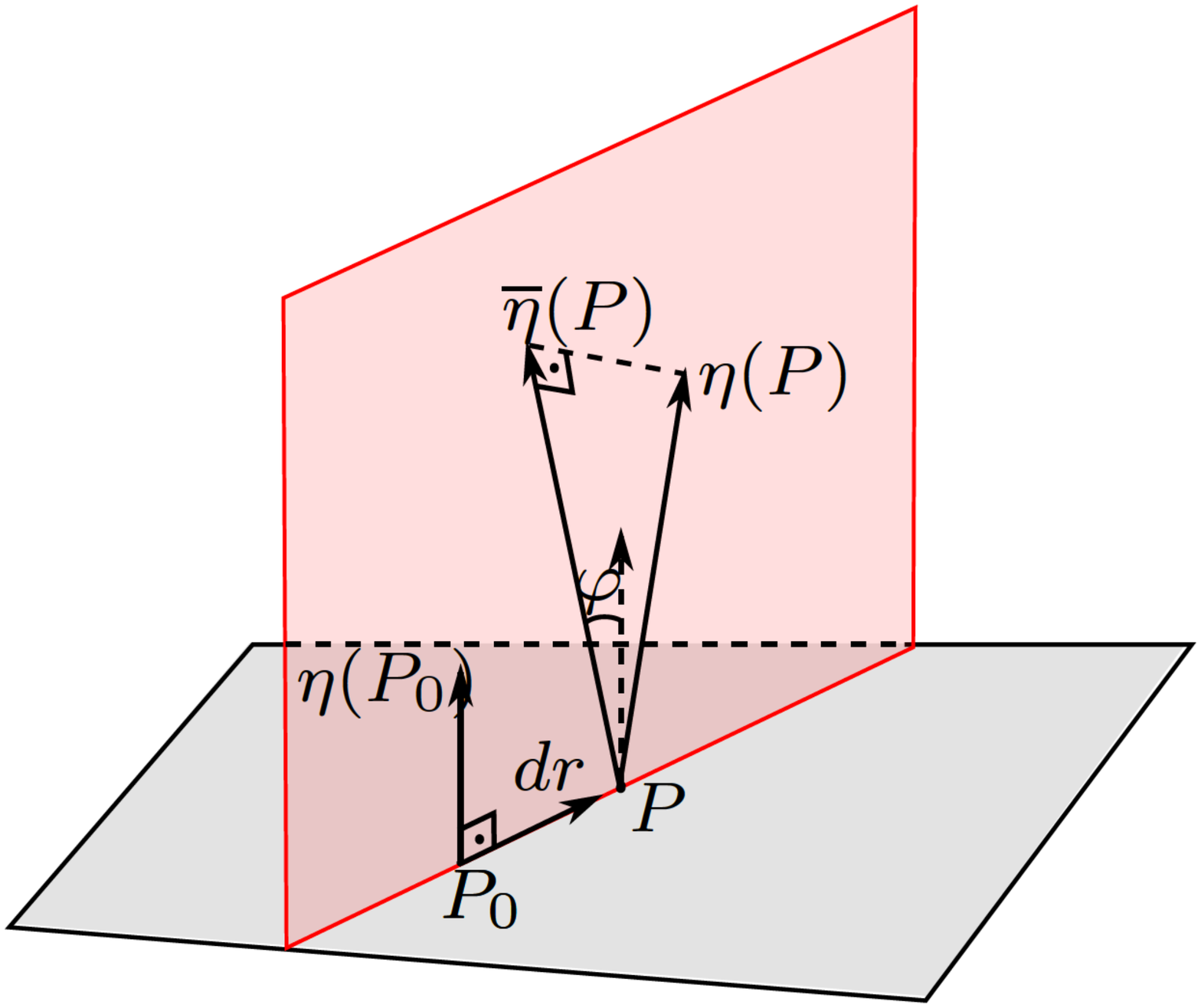}
		\caption{Geometric definition of normal curvature of the unit vector field $\eta$.}
		\label{campo}
	\end{center}
\end{figure}

Consider a plane $\Delta_0$ of the distribution $\Delta_{\eta}$ passing through a point $ P_0$ and a vector $dr=(ds,dv,dw)$, written in classical notation,    contained in the plane such that $P_0 P=P-P_0$ is in the direction of $dr$.   At the point $P$, take the orthogonal projection $\overline{\eta}(P)$ over the plane defined by the normal vector $\eta(P_0)$ and the direction $dr$ which determines an angle $\varphi$ with $\eta(P_0)$, as shown in the Fig.  \ref{campo}. For an unit vector field $\eta$, we define the  normal curvature in the direction $dr$, or normal curvature of the distribution $\Delta_{\eta}$, by
\begin{eqnarray}\label{defc}
\displaystyle k_{\eta}(dr)=\lim_{P\rightarrow P_0}{\dfrac{-\varphi}{|dr|}}.
\end{eqnarray}

When ${P\rightarrow P_0}$ we have that $\sin{\varphi} $ and $\varphi$ are of the same infinitesimal order. Once $\overline{\eta}(P)$ and the projection of $\eta(P)$ on the plane defined by $\eta(P_0)$ and the directional vector $dr$ we have that $ \dfrac{\langle \overline{\eta}(P),dr\rangle}{|\overline{\eta}(P)|}=
\dfrac{\langle \eta(P),dr\rangle}{|\eta(P)|}$ and $\sin{\varphi}= \dfrac{\langle \overline{\eta}(P), dr\rangle }{|\overline{\eta}(P)|\cdot|dr|}$, in a neighborhood of $P_0$ and so we can write 
\begin{eqnarray}\label{np}
\eta(P)=\eta(P_0)+D\eta\cdot dr+\mathcal{O}(dr),
\end{eqnarray}
where $D\eta$ denotes the differential of the vector field $\eta$ { and  $D\eta\cdot v $  denotes the derivative of $\eta$ in the direction of $v$}. It follows that 
\begin{eqnarray*}
	k_{\eta}(dr) & = & \displaystyle \lim_{P\rightarrow P_0}{\frac{-\varphi}{|dr|}}=\lim_{P\rightarrow P_0}\frac{-\sin{\varphi}}{|dr|}=\lim_{P\rightarrow P_0}
	\frac{-\langle \overline{\eta}(P),dr\rangle} {|\overline{\eta}(P)|\cdot |dr|}\cdot \frac{1}{|dr|} \\
	& =  &  \displaystyle \lim_{P\rightarrow P_0}\frac{-\langle \eta(P_0)+D\eta\cdot dr+{\mathcal O}(dr) ,dr\rangle}{|\overline{\eta}(P)|\cdot |dr|^2}=-\frac{\langle D\eta\cdot dr,dr\rangle}{\langle dr,dr\rangle},
\end{eqnarray*}
since $|\overline{\eta}(P)| \longrightarrow |\eta(P_0)|=1 $ as
$P\rightarrow P_0$. For more details see   Y. Aminov  \cite[page 8]{YA-2000}. 

The normal curvature of the vector field $\eta$ or of the plane distribution $\Delta_{\eta}$
\begin{eqnarray}\label{kn}
k_{\eta}(dr)=-\dfrac{\langle D\eta\cdot dr,dr\rangle}{\langle dr,dr\rangle}, 
\end{eqnarray}
depends only on the direction $ dr $ and not on its norm. In the holonomic case,  that is, when the plane distribution $\Delta_{\eta}$ is completely integrable, the normal curvature of the plane field in the direction $dr$ coincides with the notion of normal curvature associated to a surface. If we consider the direction $dr$ restricted to  the unit circle $  \mathbbm{S}^1$  there are two directions which corresponds to the maximum and to the minimum curvature of the plane  field and these directions can be determined as  follows. 

Given an unit   vector field  $\eta=(\eta_1,\eta_2,\eta_3)$, with $\eta_i= \eta_i(r) $, $ r=(s,v,w)\in \mathbbm{E}^3$,  $i=1,2,3$, class $ C^k$ functions, $k\geqslant 2$, denote by  $D\eta\cdot dr=$ \linebreak$(d\eta_1(r),d\eta_2(r),d\eta_3(r))$, $dr=(ds,dv,dw)$ and $d\eta_i(dr)=  \langle \nabla\eta_i,dr\rangle $. Consider the multiplier Lagrange problem:

\begin{equation}\label{eqLagrange}
\nabla k_{\eta}(dr)=\lambda\cdot \nabla G(dr), \ \ G(dr)=\langle \eta, dr \rangle.
\end{equation}
Differentiation of equation  \eqref{kn} with respect to the variable $dr$ leads to
\begin{eqnarray*}
	Dk_{\eta}\cdot dr(v) = &-\displaystyle \dfrac{\left(\langle D\eta\cdot v, dr \rangle +\langle D\eta\cdot dr, v\rangle \right)\cdot \langle dr,dr\rangle-
		2\langle dr,v\rangle\langle D\eta\cdot dr,dr\rangle}{{\langle dr,dr\rangle}^2}\\
	= & \displaystyle\dfrac{1}{{\langle dr,dr\rangle}^2}\left\langle 2\langle
	D\eta \cdot dr,dr\rangle\cdot dr-\langle dr,dr\rangle\cdot  (D\eta+D\eta^{t})\cdot dr , v \right\rangle.
\end{eqnarray*}
Thus, 
$$\nabla k_{\eta}(dr)=\dfrac{1}{{\langle dr,dr\rangle}^2}\left( 2\langle D\eta\cdot dr,dr\rangle\cdot dr-{\langle dr,dr\rangle} \cdot (D\eta+D\eta^{t})\cdot dr \right).$$
As $\nabla G(dr)=\eta$, we get the eigenvector equation

\begin{equation}\label{eq-autovetor1}
\dfrac{1}{\langle dr,dr\rangle^2}\left( 2\langle D\eta\cdot dr,dr\rangle\cdot dr- \langle dr,dr \rangle\cdot (D\eta+D\eta^{t})\cdot dr\right)=\lambda\cdot\eta .
\end{equation}
As $|\eta|=1$ and using the fact that $\langle \eta,dr \rangle=0$, computing the inner product between
equation  (\ref{eq-autovetor1}) and $\eta$ we get
\begin{eqnarray}
\lambda=-\left\langle \dfrac{(D\eta+D\eta^{t})\cdot dr}{\langle dr,dr\rangle},\eta \right\rangle,
\end{eqnarray}
which together with   equation  \eqref{eq-autovetor1} gives
{\small \begin{eqnarray}\label{eq-autovetor2}
	2\langle D\eta\cdot dr,dr\rangle\cdot dr-\langle dr,dr\rangle\cdot (D\eta+D\eta^{t})\cdot dr +\langle dr,dr\rangle \langle(D\eta+D\eta^{t})\cdot dr,\eta \rangle\cdot\eta =0.
	\end{eqnarray}}

Evaluating  the inner product of equation (\ref{eq-autovetor2}) with the vector $dr\wedge\eta$, where $\wedge$ denotes the usual wedge product in $\mathbbm{R}^3$, it is obtained another equivalent form to the equation of eigenvectors 

\begin{equation}\label{misto1}
((D\eta+D\eta^{t})\cdot dr,dr,\eta)=0.
\end{equation}
Using   that $ (D\eta^{t}-D\eta)\cdot dr=dr\wedge {\rm{curl}}(\eta)$, where $\rm{curl}(\eta)$ denotes the curl(rotational) of the vector field $\eta$ and $D\eta^{t}$ the transpose of $D\eta$, equation  \eqref{misto1}  can be written as
\begin{eqnarray*}
	(2D\eta\cdot dr+D\eta^{t}\cdot dr-D\eta\cdot dr,dr,\eta) =& 0\\
	2(D\eta\cdot dr,dr,\eta) +(dr\wedge {\rm{curl}}(\eta), dr,\eta)=& 0\\
	2(D\eta\cdot dr,dr,\eta) +\left\langle \langle dr, dr\rangle \cdot {\rm{curl}}(\eta)-\langle dr, {\rm{curl}}(\eta)\rangle \cdot dr, \eta\right\rangle =& 0,
\end{eqnarray*} and so the system
\begin{align}\label{eq-rot}
2(D\eta\cdot dr,dr,\eta) +\left\langle {\rm{curl}}(\eta), \eta\right\rangle\cdot \langle dr,dr\rangle =& 0\\
\langle dr, \eta\rangle= & 0 \nonumber
\end{align}
or equivalently, 
\begin{align}\label{misto}
((D\eta+D\eta^{t})\cdot dr,dr,\eta)=& 0\\
\langle dr, \eta\rangle= & 0 \nonumber
\end{align}
characterizes the maximum and minimum directions of the operator $k_{\eta}$ restricted to the plane defined by   equation $\langle dr, \eta\rangle= 0 $.

Note that the term $(D\eta\cdot dr,dr,\eta)$ which appears in the first equation of the system \eqref{eq-rot} is related to the geodesic torsion $\tau_{g}=\dfrac{(D\eta\cdot dr,dr,\eta)}{\langle dr,dr\rangle}$ of the vector field $\eta$ in the direction $dr$,  see \cite[page 49]{YA-2000}. Therefore, the system \eqref{eq-rot} can be rewritten as
\begin{align}\label{eq-tg}
2\tau_{g} +\left\langle {\rm{curl}}(\eta), \eta\right\rangle =& 0\\
\langle dr, \eta\rangle= & 0 \nonumber.
\end{align}

In the next lemma we will show that the vector field $\eta$ does not need to be unit  in order to study the system of equations  \eqref{eq-rot}.

\begin{lemma}\label{extensao}
	The system of equations (\ref{eq-rot}) does not change when $\eta$ is multiplied by a non null scalar function $f(r)=f(s,v,w)$. 
\end{lemma}

\proof
Multiplying the vector field $\eta$ by the scalar function $f(r)$  it follows that
{\small
	\begin{eqnarray*}
		(2 D(f(r)\eta)\cdot dr,dr,f(r)\eta) +\left\langle {\rm{curl}}(f(r)\eta), f(r)\eta\right\rangle\cdot \langle dr,dr\rangle=& 0\\
		2((\nabla f(r)\cdot dr)\eta+f(r)D\eta\cdot dr,dr,f(r)\eta) +\langle\nabla f(r)\wedge \eta+ f(r){\rm{curl}}(\eta), f(r)\eta\rangle\cdot \langle dr,dr\rangle=& 0\\
		2(f(r))^2(D\eta\cdot dr,dr,\eta) +(f(r))^2\left\langle {\rm{curl}}(\eta),\eta\right\rangle\cdot \langle dr,dr\rangle =& 0\\
		(f(r))^2\cdot \left[( 2D\eta\cdot dr,dr,\eta) +\left\langle {\rm{curl}}(\eta),\eta\right\rangle\cdot \langle dr,dr\rangle \right]=& 0.
\end{eqnarray*} } As $\langle dr, \eta\rangle=0$ and  $ f(r)\neq0$, the Lemma is proved.\hfill $ \square $

\

Recall that $  \dfrac{\langle D\eta\cdot dr,dr\rangle}{\langle dr,dr\rangle}=-k_{\eta} $ and let $ \bigg(\dfrac{D\eta+D\eta^{t}}{2}\bigg)\cdot dr=S(dr)$. We can   interpret the equation (\ref{eq-autovetor1}) as an eigenvector equation:
\begin{equation}\label{eq-autovetor3}
S(dr)-\langle S(dr),\eta\rangle  \eta =-k_{\eta}\cdot dr.
\end{equation}

The next proposition describes some properties of the operator which defines the left hand side of the eigenvector equation(\ref{eq-autovetor3}).

\

\begin{proposition}\label{operadorP} 
	Let $P :\mathbbm{E}^3\longrightarrow \mathbbm{E}^3$ be defined by $  P(v)=S(v)-\langle S(v),\eta\rangle\cdot \eta$, 
	where $S : \mathbbm{E}^3\longrightarrow \mathbbm{E}^3 $  is a linear symmetric operator and $\pi=\{v\in \mathbbm{E}^3: \langle v,\eta\rangle=0\}$.  Then  the following holds:
	\begin{itemize}
		\item [{\rm  1)}]  The kernel of the  operator $P$, which is denoted by  ker$(P)$, is generically ker$(S) + S^{-1}(\eta)$   and  in general is   transversal to the plane $\pi$.
		\item[{\rm  2)}] The plane $\pi$ is invariant by  $P$, that is, $P(\pi)\subseteq \pi$ and $P$ restricted to $\pi$ is a symmetric operator.
	\end{itemize}
	
\end{proposition}

\proof

The Item 1) follows directly from the definition of $P$. To conclude 2), we observe that for $u,v \in \pi$,    $\langle P(u),\eta \rangle=\langle P(v),\eta \rangle=0$ and $\langle P(u), v\rangle=\langle S(u), v\rangle=\langle u, S(v)\rangle =\langle u, P(v)\rangle.$
\hfill $ \square $

\

The  Lagrange multiplier problem, restricted to the plane, defined by  equation \eqref{eqLagrange}  has two critical solutions   $k_1=k_1(p) $ and $k_2=k_2(p)$ whose opposite  values are  called $\eta-$\textit{{principal values}}  or $\eta-$\textit{{principal curvatures}}.

We associate to them three directions which are called  $\eta-$\textit{{principal directions}}.  The tangent curves to these directions are called $\eta-$\textit{{principal curvature lines}}  or simply  \textit{ $\eta-$prin\-ci\-pal lines}.

The points at which the $\eta$-principal curvatures coincide are called $\eta-$\textit{{umbilic points}} and its set will be denoted by $ \mathcal{U}(\eta)$,  the points where only two principal curvatures coincide are called  $\eta-$\textit{{partially umbilic points}}. 

The set of partially umbilic points of the plane field $\Delta_{\eta}$ is defined by $\mathcal{P}(\eta)=\{p\in \mathbbm{E}^3 | \,  k_1(p)= k_2(p)\}$.

Associated to the  $\eta-$principal curvatures $k_1(p)$ and $k_2(p)$, the  $\eta-$principal  lines define two orthogonal foliations in $\mathbbm{E}^3-\mathcal{P}(\eta)$, which are denoted respectively by $\mathcal{F}_1(\eta) $ and $\mathcal{F}_2(\eta). $ They are called $\eta-$\textit{{principal foliations}}    of the plane field $\Delta_{\eta}$. The $\eta-$\textit{principal configuration} is the triple $\mathbbm{P}(\eta)= \{\mathcal{F}_1(\eta),\mathcal{F}_2(\eta), \mathcal{P}(\eta)\}.$

The foliation defined by the integral curves of the vector field $\eta$, denoted by $\mathcal{F} _{\eta}$, is orthogonal to both  $\eta-$principal foliations.

The system of equations (\ref{eq-rot}) and (\ref{misto}) also characterizes the $\eta-$principal lines associated to  the vector field $\eta$ or to the distribution $\Delta_{\eta}$.

The eigenvector equation obtained in \eqref{misto}  together with the condition  $\langle\eta,dr\rangle=0$, is equivalent to the system of implicit differential equations 
\begin{subequations}
	\begin{alignat}{2}
	L_1\cdot ds^2+L_2\cdot dsdv+L_3\cdot dsdw+L_4\cdot dv^2+L_5\cdot dvdw+L_6\cdot dw^2 & =0\label{eq-Li}\\ 
	\eta_1\cdot ds+\eta_2\cdot dv+\eta_3\cdot dw & = 0, \label{eq-Pni}
	\end{alignat} 
\end{subequations}
where, $L_i$, $ i=1,2,3,4,5,6$ are
\begin{align}\label{eq:L123456}
L_1 & =\eta_1\left(\dfrac{\partial\eta_3}{\partial v}-\dfrac{\partial\eta_2}{\partial w}\right)
+\eta_2\left(\dfrac{\partial\eta_3}{\partial s} +\dfrac{\partial\eta_1}{\partial w} \right)
-\eta_3\left(\dfrac{\partial\eta_2}{\partial s} +\dfrac{\partial\eta_1}{\partial v} \right),\nonumber \\
L_2 & =-2\eta_1\dfrac{\partial\eta_3}{\partial s}+2\eta_2\dfrac{\partial\eta_3}{\partial v}	-2\eta_3\left(\dfrac{\partial\eta_2}{\partial v} -\dfrac{\partial\eta_1}{\partial s} \right),  \nonumber\\
L_3 & =2\eta_1\dfrac{\partial\eta_2}{\partial s}+2\eta_2\left(\dfrac{\partial\eta_3}{\partial w} -\dfrac{\partial\eta_1}{\partial s} \right)-2\eta_3\dfrac{\partial\eta_2}{\partial w},\\
L_4 & =-\eta_1\left(\dfrac{\partial\eta_3}{\partial v}+\dfrac{\partial\eta_2}{\partial w}\right)
-\eta_2\left(\dfrac{\partial\eta_3}{\partial s} -\dfrac{\partial\eta_1}{\partial w} \right)
+\eta_3\left(\dfrac{\partial\eta_2}{\partial s} +\dfrac{\partial\eta_1}{\partial v} \right), \nonumber \\	
L_5 & =2\eta_1\left(\dfrac{\partial\eta_2}{\partial v} -\dfrac{\partial\eta_3}{\partial w} \right)-2\eta_2 \dfrac{\partial\eta_1}{\partial v} +2\eta_3\dfrac{\partial\eta_1}{\partial w}, \nonumber \\
L_6 & =\eta_1\left(\dfrac{\partial\eta_3}{\partial v}+\dfrac{\partial\eta_2}{\partial w}\right)
-\eta_2\left(\dfrac{\partial\eta_3}{\partial s} +\dfrac{\partial\eta_1}{\partial w} \right)
+\eta_3\left(\dfrac{\partial\eta_2}{\partial s} -\dfrac{\partial\eta_1}{\partial v} \right), \nonumber	
\end{align}

Equivalently,   assuming that $\eta_1\neq 0$, we can determine $ds$ in  (\ref{eq-Pni}), and thus obtain the system  
\begin{subequations}
	\begin{alignat}{2}
	L(s,v,w)dw^2+M(s,v,w)dwdv+N(s,v,w)dv^2=0,\label{LMNa}\\
	\eta_1(s,v,w)ds+\eta_2(s,v,w)dv+\eta_3(s,v,w)dw=0.\label{LMNb}
	\end{alignat} 
\end{subequations}
where, $L, M, N
$  are  
\begin{alignat*}{2}
L(s,v,w) & =\eta_1(\eta_1^2+\eta_3^2)\left(\dfrac{\partial\eta_3}{\partial v} +\dfrac{\partial\eta_2}{\partial w} \right)-\eta_3(\eta_1^2+\eta_3^2)\left(\dfrac{\partial\eta_1}{\partial v} +\dfrac{\partial\eta_2}{\partial s} \right)\\
& +\eta_2(\eta_3^2-\eta_1^2)\left(\dfrac{\partial\eta_1}{\partial w} +\dfrac{\partial\eta_3}{\partial s} \right)+ 2\eta_1\eta_2\eta_3\left(\dfrac{\partial\eta_1}{\partial s} -\dfrac{\partial\eta_3}{\partial w} \right),
\end{alignat*} 
\begin{alignat*}{2}
M(s,v,w) & =-2\eta_2(\eta_1^2+\eta_3^2)\left(\dfrac{\partial\eta_1}{\partial v} +\dfrac{\partial\eta_2}{\partial s} \right)+2\eta_3(\eta_1^2+\eta_2^2)\left(\dfrac{\partial\eta_3}{\partial s} +\dfrac{\partial\eta_1}{\partial w} \right)\\
& +2\eta_1(\eta_2^2-\eta_3^2)\dfrac{\partial\eta_1}{\partial s} +2\eta_1(\eta_1^2+\eta_3^2)\dfrac{\partial\eta_2}{\partial v}
- 2\eta_1(\eta_1^2+\eta_2^2)\dfrac{\partial\eta_3}{\partial w},
\end{alignat*} 
\begin{alignat*}{2}
N(s,v,w) & =-\eta_1(\eta_1^2+\eta_2^2)\left(\dfrac{\partial\eta_2}{\partial w} +\dfrac{\partial\eta_3}{\partial v} \right)+\eta_3(\eta_1^2-\eta_2^2)\left(\dfrac{\partial\eta_1}{\partial v} +\dfrac{\partial\eta_2}{\partial s} \right)\\
& +\eta_2(\eta_1^2+\eta_2^2)\left(\dfrac{\partial\eta_1}{\partial w} +\dfrac{\partial\eta_3}{\partial s} \right)+ 2\eta_1\eta_2\eta_3\left(\dfrac{\partial\eta_2}{\partial v} -\dfrac{\partial\eta_1}{\partial s} \right).
\end{alignat*}

The singular set defined by $L(s,v,w)=$ $M(s,v,w)=$ $N(s,v,w)=0$ is the set of partially umbilic points.


In \cite{AG-2016} was studied the   behavior of the $\eta-$principal configuration $\mathbbm{P}(\eta)$ near the partially umbilic set when this set is a regular curve and under generic conditions.  This configuration resemble that of principal curvature lines of hypersurfaces of $\mathbbm{E}^4$ near Darbouxian partially umbilic points, see \cite{drs-2015}.

\begin{remark}\label{rm:integrable}
When the plane distribution $\Delta_{\eta}$ is Frobenius integrable, characterized by the condition $\langle {\rm{curl}}(\eta), \eta \rangle = 0$, we are in the case of principal curvature lines of   surfaces in $\mathbbm{E}^3$.
\end{remark}


\section{First return map associated to a $\eta$-principal cycle}\label{LCfechada-ret}

In this section we will analyze the first return map associated to   a closed  $\eta$-principal line $ \gamma$, determined by the system  \eqref{misto}. 
We will assume that the curve $ \gamma\colon \mathbbm {R} \longrightarrow \mathbbm {E}^3 $ is parametrized by the arc length $ s $ and  has length $ L$.

We will describe in this section the  Poincar\'e map, or first return map, associated to   a $\eta$-principal cycle $\gamma $ of the  foliation $\mathcal{F}_1(\eta) $. Let us suppose, to fix the notation, that $\gamma$  is an orbit of a vector field $ X_1 $ defined in a tubular neighborhood of $\gamma$ and  belonging to the  plane field $\Delta_{\eta}$, that is, $X_1\in \Delta_{\eta}$.

Consider a vector field $\eta \in \mathfrak{X}_{\mathfrak{R}}^r(\mathbbm{E}^3)$, $r\geqslant 2$,   and let $\gamma$ be a closed $\eta-$principal line of the $\eta-$principal foliation $\mathcal{F}_1(\eta) $  and  $X_1$ be a vector field in the  plane distribution $\Delta_{\eta}$, such that $X_1(s)=X_1(\gamma(s))=\gamma^{\prime}(s)$, and $ X_2(s) $ be the normal unit vector along the curve $\gamma$,  $L-$periodic,  such that $\{ X_1(s),X_2(s)\}$ is a positively oriented basis of the plane of the distribution $\Delta_{\eta}$ which passes through $\gamma(s)$. Define a positively oriented orthonormal frame along $\gamma(s)$ given by $\{ X_1(s),X_2(s),N(s)\}$, where $N(s)=X_1(s)\wedge X_2(s)$ and $N(s)=N(\gamma(s))=\eta(s)$. The    Darboux equations are given by
\begin{align}
DX_1\cdot X_1 & =  k_{1}X_2+k_2N, \nonumber\\
DX_2\cdot X_1 & =  -k_{1}X_1+k_3N, \\
DN\cdot X_1 & = - k_{2}X_1-k_3X_2\nonumber. 
\end{align}

In the sequence  will assume that the frame $\{ X_1(s),X_2(s), N(s)\}$ is   $L$-periodic, where $L $ is the length of $\gamma.$ Taking a double covering, always there exists a frame $\{ X_1(s),X_2(s),N(s)\}$ which is  $2L$-periodic. In this case it is necessary to consider the second return   Poincar\'e map. See remark \ref{rem:mobius}.

Let $V_{\delta}(\gamma)$ be a tubular neighborhood of the integral curve  $\gamma$ as above and a parametrization $\alpha$ in the chart  $(s,v,w)$, $L-$periodic in the variable $s$, given by 
\begin{eqnarray}\label{parmvis}
\alpha(s,v,w)=\gamma(s)+v\cdot X_2(s)+w\cdot N(s),
\end{eqnarray}
whose Jacobian matrix $D\alpha(s,v,w)$  relative to  the basis $ \{ ds,dv,dw\} $ and $ \{ X_1, $ $ X_2, N\} $ is given by
\begin{equation*}
D\alpha(s,v,w)=
\left(
\begin{array}{ccc}
1-k_1(s)v-k_2(s)w & 0  &  0\\ 
-k_3(s)w & 1 &  0\\
k_3(s)v & 0 & 1\\
\end{array}
\right).
\end{equation*}

Let $\mathbf{X}_1$  and $\mathbf{X}_2$ be $C^r$ local vector fields generating the plane field distribution $\Delta_\eta$ in a neighborhood $V_{\delta}(\gamma)$ of $\gamma$.

As $p=\alpha(s,v,w)$ and $dp=D\alpha(s,v,w)\cdot (ds,dv,dw)$,
it follows from Hadamard's Lemma that $\mathbf{X}_1$  and $\mathbf{X}_2$ are given in the chart $(s,v,w)$ by:
{\footnotesize  \begin{align}
\mathbf{X}_1(p)  & =  X_1(s)\nonumber\\
& +\bigg (A_1(s)v+A_2(s)w+\dfrac{1}{2}A_{10}(s)v^2+A_{11}(s)vw+\dfrac{1}{2}A_{01}(s)w^2+{\mathcal O}(3)\bigg)X_2(s)\nonumber\\
& +  \bigg( B_1(s)v+ B_2(s)w+\dfrac{(1}{2}B_{10}(s)v^2+B_{11}(s)vw+\dfrac{1}{2}B_{01}(s)w^2+{\mathcal O}(3)\bigg) N(s), \label{parmX1}\
\end{align}
\vspace{-0,5cm}
\begin{align}
\mathbf{X}_2(p) &  =( C_1(s)v+ C_2(s)w+\dfrac{1}{2}C_{10}(s)v^2+C_{11}(s)vw+\dfrac{1}{2}C_{01}(s)w^2+{\mathcal O}(3)\bigg)X_1(s)\nonumber\\
& +\bigg(1+E_1(s)v+ E_2(s)w+\dfrac{1}{2}E_{10}(s)v^2+E_{11}(s)vw+\dfrac{1}{2}E_{01}(s)w^2+{\mathcal O}(3)\bigg)X_2(s) \nonumber\\
& +\bigg(F_1(s)v+ F_2(s)w+\dfrac{1}{2}F_{10}(s)v^2+F_{11}(s)vw+\dfrac{1}{2}F_{01}(s)w^2+{\mathcal O}(3)\bigg)N(s) \label{parmX2}.
\end{align}
}

The unit vector field $\mathbf{N}(p)=N(\alpha(s,v,w))$ in the neighborhood $V_{\delta}(\gamma)$ is therefore given by
{\small
\begin{eqnarray}\label{paramNLC}
\mathbf{N}(p)=\dfrac{\mathbf{X}_1(p)\wedge \mathbf{X}_2(p)}{|\mathbf{X}_1(p)\wedge \mathbf{X}_2(p)|}=\mathbf{N}_1(p)\cdot X_1(s)+\mathbf{N}_2(p)\cdot X_2(s)+\mathbf{N}_3(p)\cdot N(s),
\end{eqnarray}
}
with 
\begin{align}
\mathbf{N}_1(p)=& -B_1(s)v-B_2(s)w+(A_1(s)F_1(s)-\frac{B_{10}(s)}{2})v^2  \nonumber\\
+ & (A_1(s)F_2(s)+A_2(s)F_1(s)-B_{11}(s))vw+ b_{01}(s)w^2 + {\mathcal O}(3),\\
\nonumber\\
\mathbf{N}_2(p)=& -F_1(s)v-F_2(s)w+f_{10}(s)v^2+f_{11}(s)vw+f_{01}(s)w^2+ {\mathcal O}(3),\\
\nonumber\\
\mathbf{N}_3(p)=& \ 1-\frac{F_1(s)^2+B_1(s)^2}{2}v^2-(F_2(s)F_1(s)+B_2(s)B_1(s))vw\nonumber\\
-& \frac{F_2(s)^2+B_2(s)^2}{2}w^2+ {\mathcal O}(3),
\end{align}
and
\begin{align*}
b_{01} = & \frac{1}{2}\big(2A_2(s)F_2(s)-B_{01}(s)\big)\\
f_{01} = & \frac{1}{2}\big(2F_2(s)E_2(s)+2B_2(s)C_2(s)-F_{01}(s)\big)\\
f_{11} = & F_1(s)E_2(s)+F_2(s)E_1(s)+B_1(s)C_2(s)+B_2(s)C_1(s)-F_{11}(s)\\
f_{10} = & \frac{1}{2}\big(2 F_1(s) E_1(s)+2 B_1(s) C_1(s)-F_{10}(s)\big).
\end{align*}
Evaluating the derivatives of $\mathbf{N}$ with respect to   $s$, $v$ and $w$, it follows that:

{\footnotesize
\begin{align*}
\dfrac{\partial}{\partial s}\mathbf{N}(p) &= \bigg(-k_2(s)+\Big(F_1(s)k_1(s)-\dfrac{d}{ds}(B_1(s)\Big)v+\Big(F_2(s)k_1(s)-\dfrac{d}{ds}(B_2(s)\Big)w\\
& +{\mathcal O}(2)\bigg)\cdot X_1(s) + \bigg(  -k_3(s)-\Big(B_1(s)k_1(s)+\dfrac{d}{ds}F_1(s) \Big)v-\Big(k_1(s)B_2(s)\\
& +\dfrac{d}{ds}(F_2(s)\Big)w  + {\mathcal O}(2) \bigg)\cdot X_2(s) 
+\bigg(\Big(-F_1(s)k_3(s)-B_1(s)k_2(s))\Big)v \\
& +\Big(-F_2(s)k_3(s)-B_2(s)k_2(s)\Big)w 
+ {\mathcal O}(2)\bigg)\cdot N(s),
\end{align*}
}
{\footnotesize
\begin{align*}
\dfrac{\partial}{\partial v}\mathbf{N}(p)=& \bigg(-B_1(s)+\Big(2A_1(s)F_1(s)-B_{10}(s)\Big)v
+\Big(A_2(s)F_1(s)+ A_1(s)F_2(s)-B_{11}(s)\Big)w\\
& + {\mathcal O}(2)\bigg)\cdot X_1(s) +\bigg(-F_1(s)+2f_{10}(s)v+f_{11}(s)w + {\mathcal O}(2)\bigg)\cdot X_2(s)\\ 
& +\bigg(-\Big(F_1(s)^2+B_1(s)^2\Big)v -\Big(F_2(s)F_1(s)+B_2(s)B_1(s)\Big)w +{\mathcal O}(2)\bigg)\cdot N(s),
\end{align*}
}
{\footnotesize
\begin{align*}
\dfrac{\partial}{\partial w}\mathbf{N}(p) = & \bigg(-B_2(s)+\Big(A_2(s)F_1(s)+A_1(s)F_2(s)
-B_{11}(s)\Big)v+2b_{01}(s)w+{\mathcal O}(2) \bigg)\cdot X_1(s)\\
+ &\bigg(-F_2(s)+ f_{11}(s)v+2f_{01}(s)w +{\mathcal O}(2)\Bigg)\cdot X_2(s)\\
+ &\bigg(-\Big(F_2(s)F_1(s)+B_2(s)B_1(s)\Big)v -\Big(F_2(s)^2+B_2(s)^2\Big)w+{\mathcal O}(2)\bigg)\cdot N(s).
\end{align*}
}

As $ p =\alpha(s,v,w)$, we have that  $dp=D\alpha(s,v,w)\cdot (ds,dv,dw)$. Evaluating the mixed product and the equation of the plane \eqref{misto} which characterizes the $\eta-$principal lines, given by $\Big( (D\mathbf{N}+D\mathbf{N}^{t})\cdot dp,dp,\mathbf{N}(p) \Big)=0$ and $\langle \mathbf{N}(p),dp\rangle=0$,   we obtain
\begin{align}
L_1\cdot ds^2+L_2\cdot dsdv+L_3\cdot dsdw+L_4\cdot dv^2+L_5\cdot dvdw+L_6\cdot dw^2&=0,\label{sis-viz-tubular}\\
M_1\cdot ds+M_2\cdot dv+M_3\cdot dw&=0\label{plano-tubular}. 
\end{align}
where $L_i=L_i(s,v,w)$, $i=1,2,3,4,5,6$ and  $M_i=M_i(s,v,w)$, $i=1,2,3$ are given by
\begin{align*}
M_1(s,v,w)&= (k_3(s)-B_1(s))v-B_2(s)w\\
& +\Big(B_1(s)k_1(s)+A_1(s)F_1(s)-\frac{B_{10}(s)}{2}\Big)v^2\\
& +(B_2(s)k_1(s)+B_1(s)k_2(s)+F_1(s)k_3(s)+A_1(s)F_2(s)\\
& +A_2(s)F_1(s) -B_{11}(s))vw\\
& +(k_2(s)B_2(s)+F_2(s)k_3(s)+b_{01}(s))w^2+{\mathcal O}(3),
\end{align*}
\begin{align*}
M_2(s,v,w)&= -F_1(s)v-F_2(s)w+f_{10}(s)v^2+f_{11}(s)vw+f_{01}(s)w^2+ {\mathcal O}(3),
\end{align*}
\begin{align*}
M_3(s,v,w)& = 1-\frac{F_1(s)^2+B_1(s)^2}{2}v^2-(F_2(s)F_1(s)+B_2(s)B_1(s))vw \\
& -\frac{F_2(s)^2+B_2(s)^2}{2}w^2+ {\mathcal O}(3).
\end{align*}

{\footnotesize
\begin{align*}
L_1(s,v,w) & = B_1(s)+k_3(s)-\Big(\big(B_1(s)+2k_3(s)\big)k_1(s)-F_2(s)k_3(s)-\big(2A_1(s)-B_2(s)\big)F_1(s)\\
& +B_{10}(s) +\dfrac{d}{ds}F_1(s)\Big)v-\Big(\big(A_2(s)-2k_3(s)\big)F_1(s)
+ (A_1(s)-B_2(s))F_2(s)\\
& +2B_1(s)k_2(s)-B_2(s)k_1(s)-B_{11}(s)-\dfrac{d}{ds}F_2(s)\Big)w+{\mathcal O}(2),
\end{align*}
\begin{align*}
L_2(s,v,w) & = +2\big(F_1(s)-k_2(s)\big)+\Big(2k_1(s)k_2(s)-\big(B_1(s)+k_3(s)\big)B_2(s)
+F_1(s)F_2(s)\\
& -4f_{10}(s)-2\dfrac{d}{ds}B_1(s)\Big)v+\Big(2B_1(s)k_3(s)+2k_2(s)^2-2k_3(s)^2-2F_1(s)k_2(s)\\
& +2F_2(s)k_1(s)+F_2(s)^2-B_2(s)^2+2f_{11}(s)-2\dfrac{d}{ds}B_2(s)\Big)w+{\mathcal O}(2),
\end{align*}
\begin{align*}
L_3(s,v,w) & = F_2(s)+\Big(2B_1(s)^2+F1(s)^2+B_1(s)k_3(s)-2F_1(s)k_2(s)-F_2(s)k_1(s)\\
& -f_{11}(s)\Big)v+\Big(F_2(s)F_1(s)+2B_1(s)B_2(s)-3F_2(s)k_2(s)+2B_2(s)k_3(s)\\
& -2f_{01}(s)\Big)w+{\mathcal O}(2),
\end{align*}
\begin{align*}
L_4(s,v,w) & = -B_1(s)-k_3(s)+\Big(2A_1(s)F_1(s)-F_2(s)B_1(s)-B_1(s)k_1(s)-B_{10}(s)\\
& -\dfrac{d}{ds}F_1(s)\Big)v+\Big(A_1(s)F_2(s)-F_2(s)B_2(s)-k_1(s)B_2(s)+A_2(s)F_1(s)\\
& -B_{11}(s)-\dfrac{d}{ds}F_2(s)\Big)w+{\mathcal O}(2),
\end{align*}
\begin{align*}
L_5(s,v,w) & = -B_2(s)+\Big(A_1(s)F_2(s)+F_1(s)B_1(s)-B_1(s)k_2(s)-2F_1(s)k_3(s)\\
& +A_2(s)F_1(s)-B_{11}(s)\Big)v+\Big(2F_1(s)B_2(s)-F_2(s)B_1(s)-2F_2(s)k_3(s)\\
& -B_2(s)k_2(s)+2b_{01}(s)\Big)w+{\mathcal O}(2),
\end{align*}
\begin{align*}
L_6(s,v,w) & = \Big(F_2(s)B_1(s)-F_1(s)B_2(s)\Big)v+{\mathcal O}(2).
\end{align*}
}

In the chart $(s,v,w)$ consider two  transversal sections, $\Sigma_1=\{s=0\}$ and   $\Sigma_2=\{s=L\}$. By construction, $\alpha(\Sigma_1)=\alpha(\Sigma_2)=\Sigma$ is a transversal section.   Let  $\gamma$ and  a tubular neighborhood $V_{\delta}(\gamma)$, being $\gamma$ a $\eta-$principal cycle    defined implicitly by the system of equations
\begin{align}
\Big( (D\mathbf{N}(p)+D\mathbf{N}^{t}(p))\cdot dp,dp,\mathbf{N}(p) \Big)& =0,\\
\langle \mathbf{N}(p),dp\rangle & = 0.
\end{align}

Define the Poincar\'e first return map in the chart $(s,v,w) $ by $\pi \colon\Sigma_1\longrightarrow \Sigma_2$, by $\pi(v_0,w_0)=\Big(v(L,v_0,w_0),w(L,v_0,w_0)\Big)$, with $v(0,v_0,w_0)=v_0$ and $w(0,v_0,w_0)=w_0$. 
In order to calculate the derivative of the Poincar\'e map, we consider the system defined by equations \eqref{sis-viz-tubular} and \eqref{plano-tubular} rewritten as

\begin{align}
L_1 + L_2\cdot \dfrac{dv}{ds}+L_3\cdot \dfrac{dw}{ds}+L_4\cdot \bigg(\dfrac{dv}{ds}\bigg)^2+L_5\cdot \dfrac{dv}{ds}\dfrac{dw}{ds}+L_6\cdot \bigg(\dfrac{dw}{ds}\bigg)^2&=0,\label{sis-derivado}\\
M_1+M_2\cdot \dfrac{dv}{ds}+M_3\cdot \dfrac{dw}{ds}&=0\label{plano-derivado}.
\end{align}

Differentiating   implicitly the equations  \eqref{sis-derivado} and \eqref{plano-derivado} in relation to the initial conditions $v_0$ and $w_0$ and  evaluating at $(s,0,0)$, we have
\begin{equation}\label{eq-variacional-lc}
\left( \begin{array}{cc}
L_2 & L_3 \\ 
\cr
M_2 &  M_3 
\end{array}
\right)\cdot
\dfrac{d}{ds}
\left(
\begin{array}{cc}
\frac{\partial v}{\partial v_0}  &  \frac{\partial v}{\partial w_0} \\ 
\cr
\frac{\partial w}{\partial v_0}  &  \frac{\partial w}{\partial w_0} 
\end{array}
\right)=
\left(
\begin{array}{cc}
-\dfrac{\partial L_1}{\partial v} &  -\dfrac{\partial L_1}{\partial w} \\ 
\cr
-\dfrac{\partial M_1}{\partial v}& -\dfrac{\partial M_1}{\partial w}
\end{array}
\right)
\left(
\begin{array}{cc}
\frac{\partial v}{\partial v_0}  &  \frac{\partial v}{\partial w_0} \\ 
\cr
\frac{\partial w}{\partial v_0}  &  \frac{\partial w}{\partial w_0} 
\end{array}
\right).
\end{equation} 

For simplicity, we write the equation \eqref{eq-variacional-lc}  as
\begin{equation}\label{eq-v-simples}
A(s)\cdot \dfrac{d}{ds}U(s)=B(s)\cdot U(s).
\end{equation}

\begin{proposition}\label{prop-lc}
Let $\gamma$ be a $\eta-$principal line parametrized by arc length $s$ of the foliation $\mathcal{F}_1(\eta).$  
Suppose that there is a neighborhood $V_{\delta}(\gamma)$ of $\gamma$ such that the orthonormal frame $\{ {X}_1,{X}_2,N\}$,  $L-$periodic, is defined along $\gamma$  and the vector fields $\mathbf{X}_1(p)$ and $\mathbf{X}_2(p)$ given by (\ref{parmX1}) and (\ref{parmX2}) are defined in a  neighborhood $V_{\delta}(\gamma)$.
Then the following holds.

\item{ \rm i) } $  B_1(s)=-k_3(s) $ and $ F_1(s)-k_2(s) \neq 0 $.

\item{ \rm ii) } The  Poincar\'e return map $\pi\colon \Sigma_1\longrightarrow \Sigma_2$,  where $\Sigma_1=\{s=0\}$ and  $\Sigma_2=\{s=L\}$ are   transversal sections to $\gamma$, in   neighborhood   $V_{\delta}(\gamma)$ is such that  $D\pi(0)=U(L)$, where $ U $ is the solution of the linear differential equation
\begin{align}\label{sist-var-LC}
\dot{U}=MU,\\
U(0)=I_2, \nonumber
\end{align}  
with
\begin{equation}\label{mALC}
M(s)=
\left(
\begin{array}{cc}
M_{11}(s)  &  M_{12}(s)\\ 
& \\
-2k_3(s)  &  B_2(s)
\end{array}
\right),
\end{equation}
and
{  \footnotesize 
	\begin{align*}
	M_{11}(s) & =\dfrac{-\Big(F_2(s)+k_1(s)\Big)k_3(s)+\Big(B_2(s)-2A_1(s)\Big)F_1(s)+B_{10}(s)+\dfrac{d}{ds}F_1(s)}{2\Big(k_2(s)-F_1(s)\Big)},
	\nonumber\\
	M_{12}(s)  &  =k_3(s)+\dfrac{k_1(s)B_2(s)-A_2(s) F_1(s)+\Big(2B_2(s)-A_1(s)\Big)F_2(s)+B_{11}(s)+\dfrac{d}{ds}F_2(s)}{2\Big(k_2(s)-F_1(s)\Big)}.
	\end{align*}
}
\end{proposition}

\proof
As $\gamma$ is an integral curve of a vector field, implicitly defined by the system of equations
\begin{align}
\Big( (D\mathbf{N}(p)+D\mathbf{N}^{t}(p))\cdot dp,dp,\mathbf{N}(p) \Big)& =0,\\
\langle \mathbf{N}(p),dp\rangle & = 0,
\end{align}
we have that $dp=(ds,0,0) $  satisfies the equation (\ref{sis-viz-tubular}), and therefore  $L_1(s,0,0)$ $=L_4(s,0,0)=0$. As    $\gamma$ is disjoint of the  partially umbilic set  we have $L_2(s,0,0)$ $\neq 0$. Therefore, it follows that $ B_1(s)=-k_3(s)  $ and $ F_1(s)-k_2(s) \neq 0, $  which proofs $  { \rm i)}$. 

To prove  item $ { \rm ii)}$, consider the system of differential equations obtained in (\ref{eq-variacional-lc}), with a  simplified notation of (\ref{eq-v-simples}). 
By item $	 { \rm i)}$ we have  that $ \det (A(s))=2\Big(F_1(s)-k_2(s)\Big)  \neq 0$ and so the matrix $A(s)$ is invertible. With the initial conditions $v(0,v_0,w_0)=v_0$ and $w(0,v_0,w_0)=w_0$ we have $ U(0)=I_2 $ and  therefore we can conclude that
\begin{align*}\label{eq-v-final}
\dfrac{d}{ds}U(s) & =(A(s))^{-1}B(s)\cdot U(s)\\
U(0)& = I_2.
\end{align*}
This leads to  the result since $ (A(s))^{-1}B(s)=M(s)$.
\hfill $ \square$

\begin{proposition}
In the case where the distribution $\Delta_{\eta}$ is completely integrable, we have:
\item { \rm i) }  $k_3(s)\equiv 0$.
\item { \rm ii) } 
The derivative of the Poincar\'e map $D\pi (0)=U(L)$  is given by 
\begin{equation}
U(L)=	\left(
\begin{array}{cc}
\frac{\partial v}{\partial v_0}(L)  &  \frac{\partial v}{\partial w_0}(L) \\ 
\frac{\partial w}{\partial v_0}(L)  &  \frac{\partial w}{\partial w_0}(L) 
\end{array},
\right)
\end{equation}
where, 	
{\footnotesize
	\begin{align}
	\frac{\partial w}{\partial v_0} (L)&=0,\nonumber\\
	\frac{\partial w}{\partial w_0}(L)& =\exp\bigg(\displaystyle\int_0^L B_2(s)ds\bigg),\nonumber\\
	\frac{\partial v}{\partial v_0}(L) & =\exp\bigg(\displaystyle\int_0^L \frac{-\frac{d}{ds}F_1(s)}{F_1(s)-k_2(s)} ds\bigg), \label{dvdvo-soto}
	\end{align}
	\begin{align*}
	\frac{\partial v}{\partial w_0}(L) & =\frac{\partial v}{\partial v_0}(L)
	\int_0^L \frac{\partial v}{\partial v_0}(s)\frac{\partial w}{\partial w_0}(s) \frac{-B_2(s)\Big(2F_2(s)+k_1(s)\Big)-\frac{d}{ds}F_2(s)}{F_1(s)-k_2(s)}\  ds,
	\end{align*}
}
with \ \ 
$\frac{\partial w}{\partial w_0}(s) =\exp\bigg(\displaystyle\int_0^s B_2(t)dt\bigg)$ and 
\ \ $\frac{\partial v}{\partial v_0}(s) =\exp\bigg(\displaystyle\int_0^s \frac{-\frac{d}{dt}F_1(t)}{F_1(t)-k_2(t)} dt\bigg)$.

\end{proposition}
\proof
From  the  equation $\omega=\langle \mathbf{N}(p),dp\rangle=0$ which defines the plane distribution  $\Delta_{\eta }$ in a tubular neighborhood of the $\eta-$principal cycle $\gamma$ and evaluating $\omega $ in the  system of coordinates $(s,v,w)$, we have
\begin{equation}\label{omega}
\omega(s,v,w)=M_1(s,v,w)\cdot ds+M_2(s,v,w)\cdot dv+M_3(s,v,w)\cdot dw
\end{equation} 
with $M_1(s,v,w)$, $M_2(s,v,w)$ and $M_3(s,v,w)$ as in equation (\ref{plano-tubular}).

Differentiating the differential form $\omega$ given by equation (\ref{omega}), we have 
\begin{equation}\label{domega}
d\omega=-\dfrac{\partial M_1 }{\partial v} ds\wedge dv-\dfrac{\partial M_2}{\partial w} ds\wedge dw+\Big(\dfrac{\partial M_3}{\partial v}-\dfrac{\partial M_2}{\partial w}\Big)dv\wedge dw.
\end{equation} 
Performing the calculations we have
\begin{align*}
\omega \wedge d\omega=f(s,v,w)\cdot ds\wedge dv\wedge dw=0,
\end{align*}
where, 
\begin{equation}
f(s,v,w)=M_1\cdot\bigg(\dfrac{\partial M_3 }{\partial v}-\dfrac{\partial M_2 }{\partial w}\bigg)+M_2\cdot\bigg(\dfrac{\partial M_1 }{\partial w}-\dfrac{\partial M_3 }{\partial s}\bigg)+M_3\cdot\bigg(\dfrac{\partial M_2 }{\partial s}-\dfrac{\partial M_1 }{\partial v}\bigg).
\end{equation}
Making use of the   condition of integrability ($\omega \wedge d\omega=0$)  we have  $f(s,v,w)\equiv 0$.  Thus, as $f(s,0,0)=-2k_3(s)=0$,   the item i) is proved. Differentiating $f$ in relation to the variables $v$ and $w$ and  evaluating at $(v,w)=(0,0)$, with $k_3(s)=0$, we have
\begin{align*}
f_v(s,0,0) & =F_1(s)\Big(B_2(s)-2A_1(s)\Big)+B_{10}(s)-\dfrac{d}{ds}F_1(s),\\
f_w(s,0,0) & =-F_2(s)A_1(s)-F_1(s)A_2(s)-B_2(s)k_1(s)+B_{11}(s)-\dfrac{d}{ds}F_2(s).
\end{align*}
Solving equations   $ f_v(s,0,0)=0 $ and $ f_w(s,0,0)=0 $, respectively in  $B_{10}(s)$ and $ B_{11}(s) $ and replacing in the system (\ref{sist-var-LC}),  with the condition $ k_3(s)=0$, it is obtained
{\small
\begin{equation*}
\dfrac{d}{ds}
\left(\begin{array}{cc}
\frac{\partial v}{\partial v_0}  &  \frac{\partial v}{\partial w_0} \\ 
\frac{\partial w}{\partial v_0}  &  \frac{\partial w}{\partial w_0} 
\end{array}
\right)=
\left(
\begin{array}{cc}
-\frac{\frac{d}{ds}F_1(s)}{F_1(s)-k_2(s)} &  \frac{-B_2(s)\Big(F_2(s)+k_1(s)\Big)-\frac{d}{ds}F_2(s)}{F_1(s)-k_2(s)} \\ 
0 & B_2(s)
\end{array}
\right)
\left(
\begin{array}{cc}
\frac{\partial v}{\partial v_0}  &  \frac{\partial v}{\partial w_0} \\ 
\frac{\partial w}{\partial v_0}  &  \frac{\partial w}{\partial w_0} 
\end{array}
\right)
\end{equation*} 
}
which can be solved using the initial conditions $ \frac{\partial v}{\partial v_0}(0,0,0)=1$, $ \frac{\partial v}{\partial w_0}(0,0,0)=0  $, $ \frac{\partial w}{\partial v_0}(0,0,0)=0$ and $ \frac{\partial w}{\partial w_0}(0,0,0)=1$, getting the result.\hfill $ \square $

\begin{remark}
It is worth to note that equation (\ref{dvdvo-soto}) which determines $ \dfrac{\partial v}{\partial v_0}(L) $ is exactly  the derivative of the Poincar\'e map associated to  a principal cycle  $c:[0,L]\rightarrow  \mathbbm{M}^2$ of principal curvature lines of surfaces obtained by     Gutierrez and Sotomayor \cite{CGJS-1982}.\end{remark}

\subsection{Hyperbolicity of a $\eta$-principal cycle}\label{sec:hip}

In this subsection we will present results on the hyperbolicity  of $\eta$-principal cycles of a plane distribution in   $\mathbb{E}^3$. We will say that a closed  $\eta$-principal line $\gamma$ of length $L$ it is a  {\textit{hyperbolic  $\eta$-principal cycle,}} if the derivative of the Poincar\'e map  $D\pi(0)=U(L)$, obtained in the proposition  \ref{prop-lc}, has no eigenvalues in the unit  circle $\mathbbm{S}^1.  $

\

\begin{lemma}{\label{per-cont}}
With the same hypotheses as in Proposition (\ref{prop-lc}), consider the  perturbation  $\eta_{\varepsilon}$ of the vector field $\eta$ in the neighborhood $V_{\delta}(\gamma)$ given by
\begin{equation}
\eta_{\varepsilon}(p)=\mathbf{N}_{\varepsilon}(p)=\frac{\mathbf{X}_{1\varepsilon}(p) \wedge \mathbf{X}_2(p)}{|\mathbf{X}_{1\varepsilon}(p) \wedge \mathbf{X}_2(p)|},
\end{equation}
with
\begin{align}
\mathbf{X}_{1\varepsilon}(p) = &    \mathbf{X}_{1}(p)+ \varepsilon\Big( \phi_1(s)\cdot w + \phi_2(s)\cdot v^2+ \phi_3(s)\cdot vw\Big)N(s) \label{paramX1e},
\end{align}
and  $\phi_ i(s)$ of class $C^\infty$  for   $i=1,2,3$. Then the conditions  $ B_1(s)=-k_3(s)  $ and  $ F_1(s)-k_2(s) \neq 0 $ are invariant by this perturbation  and $\gamma$ is a closed  $\eta$-principal line of length $L$ of $\mathcal{F}_1(\eta_{\varepsilon}) $. The derivative of the Poincar\'e map $\pi_{\varepsilon}$ defined in $V_{\delta}(\gamma) \cap \Sigma$, where $ \Sigma $ is a   transversal section to $\gamma$, is given by  $D\pi_{\varepsilon}(0)=V(L)$, and  $V$ is a solution to Cauchy problem 
\begin{align}{\label{prob-cV}}
\dot{V}(s)= & M_{\varepsilon}(s)V (s)\nonumber,\\ 
V(0)= & I_2,  \ \ M_{\varepsilon}(s)=M_{\varepsilon}(s+L),
\end{align}  
with
{\small
	\begin{equation}\label{mAe}
	M_{\varepsilon}(s)=M(s)+
	\varepsilon\cdot \left(
	\begin{array}{cc}
	\dfrac{\phi_1(s) F_1(s)+2\phi_2(s) } {2\Big(k_2(s)-F_1(s)\Big)} &  \dfrac{\Big(2 F_2(s)+k_1(s)\Big)\phi_1(s)+\phi_3(s) } {2\Big(k_2(s)-F_1(s)\Big)} \\ 
	& \\
	0 &   \phi_1(s)
	\end{array}
	\right),
	\end{equation} }
and 	$  M(s) $ is given by   equation \eqref{mALC}  of Proposition  \ref{prop-lc}.
\end{lemma}
\proof
Consider the vector fields $\mathbf{X}_{1\varepsilon}(p)$ (\ref{paramX1e}) and  $\mathbf{X}_{2}(p)$ (\ref{parmX2}) and   the vector field
\begin{align*}
\mathbf{N}_{\varepsilon}(p)= &  \mathbf{N}(p)+\varepsilon\bigg(  -\phi_1(s)w-\phi_2(s)v^2-\phi_3(s)vw+{\mathcal O}(3)\bigg)\cdot X_1(s)\\
& +\varepsilon\bigg(\phi_1(s)C_1(s)vw+\phi_1(s)C_2(s)w^2+{\mathcal O}(3)\bigg)\cdot X_2(s)\\
& \varepsilon\bigg(-\phi_1(s)B_1(s)vw-\Big( \phi_1(s)B_2(s)-\frac{\phi_1(s)^2\varepsilon}{2}\Big)w^2\bigg)\cdot N(s).
\end{align*}

Calculating the derivatives of $\mathbf{N}_{\varepsilon}(p)$ in the directions $s$, $v$ and $w$, respectively, we have
{\small \begin{align*}
\dfrac{\partial}{\partial s}\mathbf{N}_{\varepsilon}(p)= & \dfrac{\partial}{\partial s} \mathbf{N}(p)+ \varepsilon\cdot \Big(-\dfrac{d}{ds}\phi_1(s)w+{\mathcal O}(2)\Big)\cdot X_1(s)- \varepsilon\cdot \Big(\phi_1(s)k_1(s)w+{\mathcal O}(2)\Big)\cdot X_2(s)\\
- & \varepsilon\cdot\Big( k_2(s)\phi_1(s)w+{\mathcal O}(2)\Big)\cdot N(s),
\end{align*} }
{\small\begin{align*}
\dfrac{\partial}{\partial v}\mathbf{N}_{\varepsilon}(p)= & \dfrac{\partial}{\partial v} \mathbf{N}(p)-\varepsilon\cdot \bigg(2\phi_2(s)v+\phi_3(s)w+{\mathcal O}(2)\bigg)\cdot X_1(s)\\
+ &\varepsilon\cdot \Big(\phi_1(s)C_1(s) w+{\mathcal O}(2)\Big)\cdot X_2(s)-\varepsilon\cdot\Big( \phi_1(s)B_1(s)w+{\mathcal O}(2)\Big)\cdot N(s)
\end{align*}}
{\small\begin{align*}
\dfrac{\partial}{\partial w}\mathbf{N}_{\varepsilon}(p) & = \dfrac{\partial}{\partial w} \mathbf{N}(p)- \varepsilon\cdot \bigg(\phi_1(s)+\phi_3(s)v
+{\mathcal O}(2)\bigg)\cdot X_1(s)\\
& + \varepsilon\cdot \Big(\phi_1(s)C_1(s) v +2\phi_1(s)C_2(s)+{\mathcal O}(2)\Big)\cdot X_2(s)\\
& -\varepsilon\bigg(\phi_1(s)B_1(s)v+\Big(2\phi_1(s)B_2(s)+\varepsilon \phi_1(s)^2\Big)w+{\mathcal O}(2)\bigg)\cdot N(s).
\end{align*}}

Calculating the mixed product and the equation of the plane that characterize the  $\eta$-principal lines obtained in equation (\ref{misto}), given, respectively, by $\Big( (D\mathbf{N}_{\varepsilon}+D\mathbf{N}_{\varepsilon}^{t})\cdot dp,dp,\mathbf{N}_{\varepsilon}(p) \Big)=0$ and $\langle \mathbf{N}_{\varepsilon}(p),dp\rangle=0$,   it follows that
\begin{align}\label{sis-viz-tubular-pert}
L_{\varepsilon1}\cdot ds^2+L_{\varepsilon2}\cdot dsdv+L_{\varepsilon3}\cdot dsdw+L_{\varepsilon4}\cdot dv^2+L_{\varepsilon5}\cdot dvdw+L_{\varepsilon6}\cdot dw^2&=0,\nonumber\\
M_{\varepsilon1}\cdot ds+M_{\varepsilon2}\cdot dv+M_{\varepsilon3}\cdot dw&=0, 
\end{align}
with $L_{\varepsilon i}=L_{\varepsilon i}(s,v,w)$, $i=1,2,3,4,5,6$ and  $M_{\varepsilon i}=M_{\varepsilon i}(s,v,w)$, $i=1,2,3$, given by
\begin{align*}
M_{\varepsilon 1}(s,v,w)&= (k_3(s)-B_1(s))v-\Big(B_2(s)+\varepsilon\phi_1(s)\Big)w+{\mathcal O}(2),\\
M_{\varepsilon 2}(s,v,w)&= -F_1(s)v-F_2(s)w+{\mathcal O}(2),\\
M_{\varepsilon 3}(s,v,w)& =1+{\mathcal O}(2),
\end{align*}
{\small\begin{align*}
L_{\varepsilon 1}(s,v,w) & = L_1(s,v,w)+\varepsilon\cdot\bigg(\Big(\phi_1(s)F_1(s)+2\phi_2(s)\Big)v+\Big(\phi_1(s)\big(k_1(s)+ F_2(s)\big)+\phi_3(s)\Big)w+{\mathcal O}(2)\bigg),
\end{align*}}
{\small\begin{align*}
L_{\varepsilon 2}(s,v,w)& =L_2(s,v,w) +\varepsilon\cdot \bigg(-\phi_1(s)\Big(B_1(s)+k_3(s)\Big)v-\Big(2\phi_1^2(s)\varepsilon+2\phi_1(s)B_2(s)\\
& +2\phi_1(s)C_1(s)+2\dfrac{d}{ds}\phi_1(s)\Big)w+{\mathcal O}(2)\bigg),
\end{align*}}
{\small\begin{align*}
L_{\varepsilon 3}(s,v,w) & = L_3(s,v,w)- \varepsilon\cdot\bigg(\phi_1(s)C_1(s)v+2\phi_1(s)\Big(C_2(s)-k_3(s)-B_1(s)\Big)w+{\mathcal O}(2)\bigg),
\end{align*}}
{\small\begin{align*}
L_{\varepsilon 4}(s,v,w) & = L_4(s,v,w)+\varepsilon\cdot\bigg(-2\phi_2(s)v-\Big(\phi_1(s)F_2(s)+\phi_1(s)k_1(s)+\phi_3(s)\Big)w+{\mathcal O}(2)\bigg),
\end{align*}}
{\small\begin{align*}
L_{\varepsilon 5}(s,v,w) & = L_5(s,v,w)-\varepsilon\cdot\bigg(\phi_1(s)+\phi_3(s)v+\phi_1(s)\Big(k_2(s)-2F_1(s)\Big)w+{\mathcal O}(2)\bigg),
\end{align*}}
{\small\begin{align*}
L_{\varepsilon 6}(s,v,w)& = L_6(s,v,w)-\varepsilon\cdot \Big(F_1(s)\phi_1(s)v+{\mathcal O}(2)\Big).
\end{align*}}

Substituting in the variational equation (\ref{eq-variacional-lc}), we get
\begin{equation}\label{eq-variacional-lc-per}
\left( \begin{array}{cc}
L_{\varepsilon2} & L_{\varepsilon3} \\ 
\cr
M_{\varepsilon2} &  M_{\varepsilon3} 
\end{array}
\right)\cdot
\dfrac{d}{ds}
\left(
\begin{array}{cc}
\frac{\partial v}{\partial v_0}  &  \frac{\partial v}{\partial w_0} \\ 
\cr
\frac{\partial w}{\partial v_0}  &  \frac{\partial w}{\partial w_0} 
\end{array}
\right)=
\left(
\begin{array}{cc}
-\dfrac{\partial L_{\varepsilon1}}{\partial v} &  -\dfrac{\partial L_{\varepsilon1}}{\partial w} \\ 
\cr
-\dfrac{\partial M_{\varepsilon1}}{\partial v}& -\dfrac{\partial M_{\varepsilon1}}{\partial w}
\end{array}
\right)
\left(
\begin{array}{cc}
\frac{\partial v}{\partial v_0}  &  \frac{\partial v}{\partial w_0} \\ 
\cr
\frac{\partial w}{\partial v_0}  &  \frac{\partial w}{\partial w_0} 
\end{array}
\right),
\end{equation}  
which we shall denote by
\begin{equation}\label{eq-v-simples-per}
A_{\varepsilon}(s)\cdot \dfrac{d}{ds}U(s)=B_{\varepsilon}(s)\cdot U(s).
\end{equation}

Evaluating in $(s,0,0)$, we have $L_{\varepsilon 1}(s) =L_2(s)  $, $ L_{\varepsilon 2}(s) =L_2(s)$, which ensure the invariance of the conditions $ B_1(s)=-k_3  $ and $ F_1(s)-k_2(s) \neq 0 $.

We also have in $(s,0,0)$ what  $ L_{\varepsilon 3}(s) =L_3(s) $, $ M_{\varepsilon 2}(s) =M_2(s)  $ and $ M_{\varepsilon 3}(s) =M_3(s)  $, and therefore the determinant $ \det (A_{\varepsilon}(s))=2\Big(F_1(s)-k_2(s)\Big)  \neq 0$.  Therefore,  the matrix $A_{\varepsilon}(s)$ is invertible. Let   $v(0,v_0,w_0)=v_0$, $w(0,v_0,w_0)=w_0$ and    $ U(0)=I_2 $. Therefore we conclude that
\begin{align*}\label{eq-v-final}
\dfrac{d}{ds}U(s) & =(A_{\varepsilon}(s))^{-1}B_{\varepsilon}(s)\cdot U(s)\\
U(0)& = I_2.
\end{align*}
This leads to  the result,  since $(A_{\varepsilon}(s))^{-1}B_{\varepsilon}(s)=M_{\varepsilon} (s)$. \hfill $ \square $

\section{Proof of Theorem \ref{th:pert-cont}}\label{sec:th1}

Suppose that $ \gamma  $ is not a   hyperbolic $\eta-$principal cycle of $\mathcal{F}_1(\eta)$ and consider a chart $(s,v,w)$, $L-$periodic in $s$, defined by equation \eqref{parmvis}.  We will show that there is an   vector field $\eta_{\varepsilon}$ sufficiently close to the vector field $\eta$, such that $\gamma$ is a hyperbolic  $\eta_{\varepsilon}-$principal cycle of $\mathcal{F}_1(\eta_{\varepsilon})$, that is, suppose that an eigenvalue of  $U(L)$  is in the unit circle $\mathbbm{S}^1$, being $U$ the solution of the Cauchy problem
\begin{align}
\dot{U}(s)=M(s)\cdot U(s),\\
U(0)=I_2, M(s)=M(s+L)\nonumber
\end{align}  
with $M(s)$ given by equation (\ref{mALC}).  The  Cauchy problem (\ref{prob-cV}) of Lemma \ref{per-cont}\ can be interpreted as a geometric control problem
{\small\begin{equation}\label{sc-lc}
	\dot{V} (s)=M(s)V(s)+ \sum_{i=1}^{3}u_i(s) E_i(s)V(s)=\left(M(s)+\sum_{i=1}^{3}u_i(s) E_i(s)\right)V(s),   \ \ k\geqslant 1,
	\end{equation}}
with controls
\begin{align*}
u_1(s) & =\varepsilon\cdot \dfrac{\phi_1(s)F_1(s)+2\phi_2(s)}{2\Big(F_1(s)-k_2(s)\Big)} , \\
u_2(s) & =\varepsilon\cdot\dfrac{\phi_1(s)\Big(2F_2(s)+k_1(s)\Big)+\phi_3(s)}{2\Big(F_1(s)-k_2(s)\Big)}, \\  
u_3(s) & =\varepsilon\cdot\phi_1(s),
\end{align*}
with $\varepsilon$ $\in$ $\mathbbm{R}$,  $\phi_i $ 
of compact support, and 
$$E_1(s)=\left(
\begin{array}{cc}
1 &  0\\ 
0& 0
\end{array} 
\right), \ \ \
E_2(s)=\left(
\begin{array}{cc}
0&  1 \\ 
0& 0
\end{array} 
\right) \ \ \ 
E_3(s)=\left(
\begin{array}{cc}
0&  0\\ 
0& 1 
\end{array} 
\right).
$$
Let  $\{B_{i}^{j}(s)\}$, $i=1,2,3$ and   $j=1,2$  defined by  
\begin{align*}
B_{i}^{1}(s)= & E_i (s),\\
B_{i}^{2}(s)= & [E_i(s), A(s)]=A(s)E_{i}(s)-E_{i}(s)A(s).
\end{align*}
So it follows that,
\begin{equation*}
B_{1}^{2}(s)=A(s)E_{1}(s)-E_{1}(s)A(s)= 
\left(
\begin{array}{cc}
0 &  *\\ 
-2k_3(s) & 0
\end{array} 
\right).
\end{equation*}
Taking  $\overline{s}$ $\in$ $[0,L]$ such that $k_{3}(\overline{s})\neq 0$,  we have that $\{E_1(\overline {s}),E_2(\overline {s}),E_2(\overline {s}),B_{1}^{2}(\overline {s})\}$   is a basis of the vector space $M_2(\mathbbm R)$(matrices of order 2 with real coefficients).  Therefore,
\begin{equation}
{\rm Span}\{B_i^j(\overline {s}) \  | \ i \in \{ 1,\ldots 3\}, j=1,2\}=T_{I_{2}}M_2(\mathbbm R),
\end{equation}
and so by the Theorem  \ref{Tcont} 
in Appendix,  the control system (\ref{sc-lc})  is controllable in $[0,L]$. Therefore,  we can get a vector field $\eta_{\varepsilon}$,  $\varepsilon-C^{r}$ close to  $\eta$, such that the eigenvalues of $D\pi_{\varepsilon}(0)=V(L)$, solution of the Cauchy problem (\ref{prob-cV})  of  Lemma  
\ref{per-cont},  are not in the unit circle $\mathbbm{S}^1$. This ends the proof.

\begin{remark}\label{rem:mobius} The case where the orthonormal frame  $\{ X_1(s),X_2(s),N(s)\}$ is   $2L$-periodic,  we can consider the second return map to develop the analysis. This case happens when the derivative of the first return map is hyperbolic  and has two negative real eigenvalues and so there are two non orientable invariant surfaces (M\"obius strips)  containing $\gamma$.  See \cite{HPS-1977} or \cite{PM-1982}.
\vspace{-0.5cm}
\begin{figure}[H]
	\begin{center}
		\includegraphics[scale=0.3]{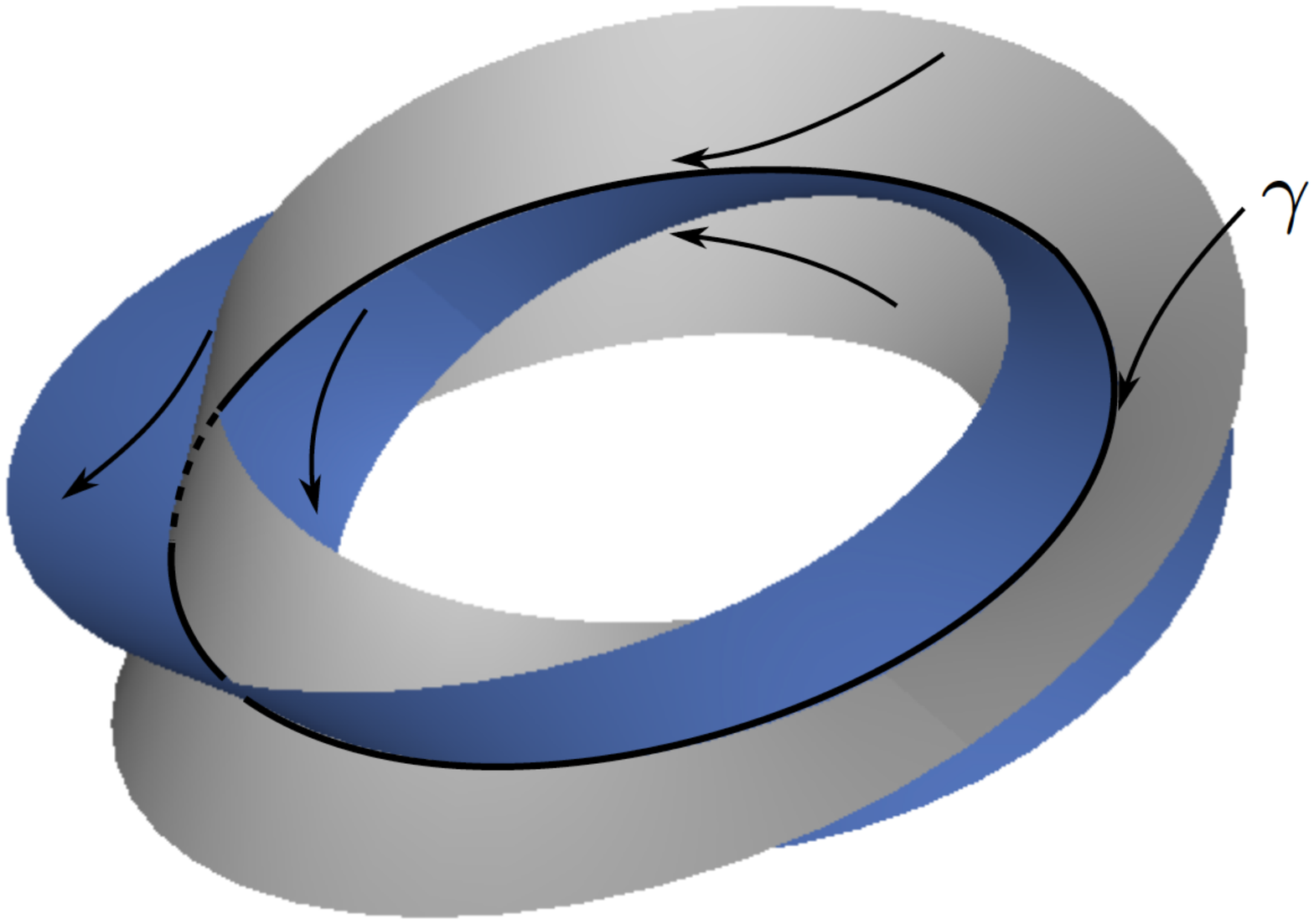}
		\caption{ Invariant surfaces having topological type of   M\"obius strips.}
		\label{Mobius}
	\end{center}
\end{figure}
\end{remark}

\begin{remark} As it is well established in geometric theory of dynamical systems, a hyperbolic compact leaf $\gamma$ of an one dimensional foliation is stable under small deformations and the invariant manifolds (stable and unstable) asymptotic to $\gamma$ are smooth, see \cite{F-1971}, \cite{HPS-1977} and \cite{PM-1982}. For higher dimensional phenomena see \cite{HM-1973} for example.
\end{remark}

\section{Proof of Theorem \ref{th:peridoc-dense}}\label{sec:th2}
Consider $ \mathbbm{E}^3$ with the exhaustion of compacts $\{ K_m, m=1,\ldots, \infty \}$ defining the $C^r$ topology, $r\geq 3.$ See \ref{app:B}.

Define the sets $\mathcal{G}(K_m,L)$ of $\eta \in\mathfrak{X}_{\mathfrak{R}}^r(\mathbbm{E}^3)$   such that all $\eta-$principal cycles  of $\mathcal{F}_i(\eta)$  ($i=1,2$), with  length smaller than $L$ and passing through the compact $K_m$ are hyperbolic.

\begin{lemma}
\label{lem:peixoto1}
The set $\mathcal{G}(K_m,L) $ is residual in $\mathfrak{X}_{\mathfrak{R}}^r(\mathbbm{E}^3)$.
\end{lemma}

\proof
It is similar to the case of vector fields, see \cite{PM-1982} and \cite{MP-1967}.\hfill $\square$

\begin{lemma}
\label{lem:peixoto2}
The set  \[\mathcal{G}= \bigcap_{L=0}^{\infty}\bigcap_{m=1}^{\infty}\mathcal{G}(K_m,L). \] is dense in $\mathfrak{X}_{\mathfrak{R}}^r(\mathbbm{E}^3)$.
\end{lemma}

\proof
This is direct consequence of the fact that $\mathfrak{X}_{\mathfrak{R}}^r(\mathbbm{E}^3)$ is a Baire space. See  \ref{app:B} or \cite{MP-1967}.\hfill $\square$

\section{Examples}\label{sec:exe}
In this section will be given examples of hyperbolic $\eta-$principal cycles of  saddle and nodal type. Also a  semi hyperbolic  $\eta-$principal cycle is given.  

Consider the vector fields $X_1=(-v+\lambda s(1-s^2-v^2),\ s+\lambda v(1-s^2-v^2) ,\ -aw)$, $X_2=(s+\varepsilon v,-\varepsilon s+v,0)\wedge X_1$ and let $\eta= (X_1\wedge X_2)/|X_1\wedge X_2|$ a unitary vector field defined in an open set of $\mathbbm{R}^3$.

\begin{proposition} In the conditions above the unit circle $\gamma(t)=(\cos t, \sin t, 0)$ is a $\eta$-principal cycle of the principal foliation $\mathcal{F}_1(\eta)$.  Moreover, the eigenvalues of the Poincar\'e map  are $\lambda_1=\exp(-4\pi\lambda)$ and  $\lambda_2=\exp(\frac{2\pi a\varepsilon(\lambda-a)}{a\varepsilon+1})$ and this cycle is:
\begin{enumerate}[i)]
	\item Hyperbolic of nodal type when $a\varepsilon>0$ and $\lambda(\lambda-a) <0$.
	\item Hyperbolic of saddle type when $a\varepsilon>0$ and $\lambda(\lambda-a)>0$.
	\item Semi hyperbolic of saddle node type when $a=0$ and  $\lambda\neq 0$ or $a\neq0$ and  $\lambda= 0$. 	
\end{enumerate}	
\end{proposition}
\proof
Consider the change of coordinates given by:  
\begin{align*}
s & = (1-v_1)\cos{t}\\
v & =  (1-v_1)\sin{t}\\
w & = w.
\end{align*}
From equations \eqref{sis-viz-tubular} and \eqref{plano-tubular}, we obtain
\begin{align}
L_1\cdot dt^2+L_2\cdot dtdv+L_3\cdot dtdw+L_4\cdot dv^2+L_5\cdot dvdw+L_6\cdot dw^2&=0,\label{mud-tubular}\\
M_1\cdot dt+M_2\cdot dv+M_3\cdot dw&=0\label{mud-plano-tubular},
\end{align}
with $L_i=L_i(t,v,w)$, $i=1,2,\ldots,6$ and  $M_i=M_i(t,v,w)$, $i=1,2,3$. Here the new variable $v_1$  is written in the old notation.
It can be verified, by straightforward calculations, that:

\[  L_3(t,0,0)=2(a\varepsilon+1), \ \  L_4(t,0,0)=-2\lambda,  \;\;\text{and}\;\; L_i(t,0,0)= 0 ,\;\; i\ne 3,4. \]
\[ M_1(t,0,0)=  M_3(t,0,0)=0  \;\; \text{and}\;\;   M_2(t,0,0)=-1,\]
\[(L_1)_v( t,0,0)=0,\; (L_1)_w(t,0,0)=2a\varepsilon (a-\lambda),\] \[(M_1)_v(t,0,0)=-2\lambda,\; (M_1)_w=(t,0,0)=0.\]
The  calculations were corroborated with symbolic computations and we omitted these long expressions. We have that
$\eta(\cos t,\sin t, 0)=(\cos t, \sin t, 0 )$ and so the   unit circle $\gamma$ is  a $\eta-$principal cycle. The $\eta-$principal curvatures are $k_1=-1$ and $k_2=a\varepsilon. $

Differentiating    equations  (\ref{mud-tubular}) and  (\ref{mud-plano-tubular}) with respect to the initial conditions $v_{0}$ and $w_0$ and  evaluating at $(t,0,0)$, we obtain  the variational equation

\begin{equation}\label{eq-variacional-lc-mud}
\left( \begin{array}{cc}
0 & 2(a\varepsilon+1) \\ 
\cr
-1 & 0 
\end{array}
\right)\cdot
\dfrac{d}{dt}
\left(
\begin{array}{cc}
\frac{\partial v}{\partial v_0}  &  \frac{\partial v}{\partial w_0} \\ 
\cr
\frac{\partial w}{\partial v_0}  &  \frac{\partial w}{\partial w_0} 
\end{array}
\right)=
\left(
\begin{array}{cc}
0 &  2a\varepsilon(\lambda-a) \\ 
\cr
2\lambda& 0
\end{array}
\right)
\left(
\begin{array}{cc}
\frac{\partial v}{\partial v_0}  &  \frac{\partial v}{\partial w_0} \\ 
\cr
\frac{\partial w}{\partial v_0}  &  \frac{\partial w}{\partial w_0} 
\end{array}
\right).
\end{equation} 
or equivalently,
\begin{equation}\label{eq-var-exe}
\dfrac{d}{dt}
\left(
\begin{array}{cc}
\frac{\partial v}{\partial v_0}  &  \frac{\partial v}{\partial w_0} \\ 
\cr
\frac{\partial w}{\partial v_0}  &  \frac{\partial w}{\partial w_0} 
\end{array}
\right)=
\left(
\begin{array}{cc}
-2\lambda &  0 \\ 
\cr
0 & \frac{a\varepsilon(\lambda-a)}{a\varepsilon+1}
\end{array}
\right)
\left(
\begin{array}{cc}
\frac{\partial v}{\partial v_0}  &  \frac{\partial v}{\partial w_0} \\ 
\cr
\frac{\partial w}{\partial v_0}  &  \frac{\partial w}{\partial w_0} 
\end{array}
\right).
\end{equation} 
Solving the linear  differential equation \eqref{eq-var-exe} it follows that:  
\[
D\pi (0)
=\left(
\begin{array}{cc}
\exp(-4\pi\lambda) &  0 \\ 
\cr
0 & \exp\Big(\frac{2\pi a\varepsilon(\lambda-a)}{a\varepsilon+1}\Big)
\end{array}
\right),
\]
which concludes the proof.\hfill $ \square $

\section{Concluding Remarks}\label{sec:cr}

The extrinsic geometry of vector fields and flows  can be applied in the dynamic of mechanical systems, control systems, theoretical  elasticity and fluid flows, see   \cite{YA-2000} and \cite{Dav-1994}.  Also, in foliation theory the extrinsic geometry of leaves and flows is  present, see \cite{RoWa-2011} for a recent exposition on this subject.

In the case where  $ \Delta_\eta$ is integrable these foliations are exactly the principal curvature lines of an one parameter family of surfaces in $ \mathbbm{E}^3,  $ a classical subject of differential geometry of surfaces which was introduced by G. Monge \cite{Mo-1796}. The qualitative theory and global aspects of principal lines were initiated  by   C. Gutierrez and J. Sotomayor \cite{CGJS-1982} and is a current subject of research.

A correlated result, which corresponds to closed principal lines of three dimensional manifolds immersed in $\mathbbm E^4$ was obtained in  \cite{RG-1993}.


The  results of this work can be extended to  regular plane fields defined  in three dimensional Riemannian manifolds, following the natural concepts of curvature tensor and second fundamental forms of plane fields \cite{Re-1977}.

\section*{Acknowledgments} This work was partially supported by PRONEX/CNPq/FAPEG.

\appendix

\section{ Control Theory Concepts}\label{app:A}

We will now discuss some preliminary results of the control theory used in the previous section. 
For a broader reading in control theory we indicate Coron[Chapter 1] \cite{JMC-2007} and Rifford[Chapter 2] \cite{Rif-2014}.

Consider a homogeneous linear differential equation   control  system \linebreak $M_{n}(\mathbbm R)$ (real matrices of order$ n\geqslant 1$) of the form 
\begin{equation}\label{sc}
\dot{X} (s)=A(s)X(s)+ \sum_{i=1}^{k}u_i(s)B_i(s)X(s),   \ \ k\geqslant 1,
\end{equation}
with $X(s)$ in $ M_n(\mathbbm R) $; the controls $u_i(s)$ in   $L^{1}([0,L]:\mathbbm R^k )$, $k\geqslant1$ and $s\in [0,L]$, $L>0$. The functions  $A(s)$,  $B_i(s)$ are continuous applications in $M_n(\mathbbm R)$, with $i=1,2,\cdots,k$. Consider   the  Cauchy problem
\begin{align}\label{sc-cauchy}
\dot{X} (s) = & A(s)X(s)+ \sum_{i=1}^{k}u_i(s)B_i(s)X(s),   \ \ k\geqslant 1, \\
X(0) = & X_0\nonumber.
\end{align}
\vspace{-0.3cm}
\begin{definition}
The control system (\ref{sc}) is completely controllable  or controllable if, for each $(X_0,X_1) \in M_n(\mathbbm R) \times M_n(\mathbbm R)$, there exists a  control \linebreak $u(t)=(u_1(s),\ldots, u_k(s)  )$ in $L^{\infty}([0,L]\colon\mathbbm R^k )$, such that the solution $X\in$ \linebreak $ C^0([0,L],M_n(\mathbbm R) )$ of the Cauchy problem (\ref{sc-cauchy}) satisfies $X(L)=X_1$.
\end{definition}
The following theorem ensure us a sufficient condition for the controllability of (\ref{sc}), and is in Rifford and Ruggiero \cite{RR-2012}. See also  See also  \cite[Chapters 1,10]{JMC-2007} .

\begin{theorem}{\label{Tcont}}
Let $L>0$, and smooth maps $s \in [0,L]\rightarrow$ $ A(s)$, $ B_1(s)\ldots $  $ B_k(s)$ $\in $ $  M_n(\mathbbm R)$. Define $k$ sequences of maps 
$$\{B_1^j\},\ldots \{B_k^j\}\colon [0,L]\rightarrow T_{I_{n}}M_n(\mathbbm R),$$ by
\begin{align}\label{seq-b}
{B_i^0(s)} := &{B_i(s)}, \nonumber \\
\vdots &  \\
{B_i^j(s)} := & \frac{d}{ds}B_i^{j-1}(s) -[B_i^{j-1}(s),A(s)] \nonumber.
\end{align}
for all $s\in[0,L]$ and each  $i\in \{1,\ldots k\}$. Suppose  that there is a $\overline {s}\in [0,L]$  such that  

\begin{equation}
{\rm Span}\{B_i^j(\overline {s}) \  | \ i \in \{ 1,\ldots k\}, j\in \mathbbm{N}\}=T_{I_n}M_n(\mathbbm R).
\end{equation}
Then, the control system (\ref{sc}) is controllable in $[0,L]$.
\end{theorem}
\section{Topology of unit vector fields in $\mathbbm{E}^3 $}\label{app:B}

Here we will present, following Peixoto \cite{MP-1967}, a topology for the space of vector fields defined in non-compact manifolds.

Let $ M^n=M $ a $n-$dimensional  non-compact, the topology  $  C^r $ of Whitney $r\geqslant 1$, in space  $\mathfrak{X}^r(M)$ of vectors fields  in $M$ is defined as follows.
 
Let
 \begin{equation}\label{exaustao}
 K_1\subset K_2\subset K_3\subset\cdots \subset K_i\subset K_{i+1}\subset \cdots \subset M,
 \end{equation}
 a decomposition of $ M $ due to compact exhaustion $ K_i $, with non-empty interior $\accentset{\circ}{K}_i\subset K_{i+1}$.
If $X$is a vector field in  $\mathfrak{X}^r(M)$ and $\delta(x)>0 $ is a positive function in $M$, Let
 \begin{equation*}
 \delta_i=\inf \left\{\delta(x), \ \ \ x\in K_i-\accentset{\circ}{K}_{i-1}\right\}, \ \ \ K_0=\phi.
 \end{equation*} 
 \begin{equation}\label{aber-top-Peixoto}
 \mathcal{A}\Big(X,\delta(x)\Big)=\bigcap_{i=1}^{\infty}\{Y:d\Big(X,Y,K_i-\accentset{\circ}{K}_{i-1}\Big)<\delta_i\},
 \end{equation}
 where $d$ is usual distance $C^r$ in compacts, see \cite{HM-1976}.
 
 The set  $\mathcal{A}\Big(X,\delta(x)\Big)$ forms a base of a neighborhood of  $X$  in  the topology in $\mathfrak{X}^r(M)$.
 
 \begin{remark}
 An important observation is that this base does not depends on the decomposition of $ M$, nor on the  chosen of metric in each compact $K_i$.
 \end{remark}

 \begin{theorem}\label{Baire}
 The set  $\mathfrak{X}^r(M)$, with the topology defined above, is a Baire space.
 	
 \end{theorem}
 \proof See \cite[Chapter 2]{HM-1976} and  \cite{MP-1967}.
  \begin{corollary}\label{BaireU}
  The set  $\mathfrak{X}_{\mathfrak{R}}^r(\mathbbm{E}^3)$ of unit regular vector fields, with the induced topology of  $\mathfrak{X}^r( \mathbbm{E}^3)$, is a Baire space.
 	
  \end{corollary}
  \proof See \cite[Chapter 2]{HM-1976}.
   \hfill $ \square $


\begin{thebibliography}{99}
	\expandafter\ifx\csname natexlab\endcsname\relax\def\natexlab#1{#1}\fi
	\expandafter\ifx\csname url\endcsname\relax
	\def\url#1{\texttt{#1}}\fi
	\expandafter\ifx\csname urlprefix\endcsname\relax\def\urlprefix{URL }\fi
	
	\bibitem[{Aminov(2000)}]{YA-2000}
	Aminov, Y., 2000. The geometry of vector fields. Gordon and Breach Publishers,
	Amsterdam.
	
	\bibitem[{Boersma and Dray(1995)}]{bo-1995}
	Boersma, S., Dray, T., 1995. Parametric manifolds. {I}. {E}xtrinsic approach.
	J. Math. Phys. 36~(3), 1378--1393.
	\newline\urlprefix\url{https://doi.org/10.1063/1.531127}
	
	\bibitem[{Coron(2007)}]{JMC-2007}
	Coron, J.-M., 2007. Control and nonlinearity. Vol. 136 of Mathematical Surveys
	and Monographs. American Mathematical Society, Providence, RI.
	
	\bibitem[{Davydov(1994)}]{Dav-1994}
	Davydov, A.~A., 1994. Qualitative theory of control systems. Vol. 141 of
	Translations of Mathematical Monographs. American Mathematical Society,
	Providence, RI, translated from the Russian manuscript by V. M. Volosov.
	
	\bibitem[{Fenichel(1971/1972)}]{F-1971}
	Fenichel, N., 1971/1972. Persistence and smoothness of invariant manifolds for
	flows. Indiana Univ. Math. J. 21, 193--226.
	
	\bibitem[{Garcia and Sotomayor(2009)}]{RGJS-2009}
	Garcia, R., Sotomayor, J., 2009. Differential equations of classical geometry,
	a qualitative theory. Publica\c c\~oes Matem\'aticas. IMPA, Rio de Janeiro,
	27${^{\rm{o}}}$ Col{\'o}quio Brasileiro de Matem{\'a}tica.
	
	\bibitem[{Garcia(1993)}]{RG-1993}
	Garcia, R.~A., 1993. Hyperbolic principal cycles on hypersurfaces of {${\bf
			R}^4$}. Ann. Global Anal. Geom. 11~(2), 185--196.
	\newline\urlprefix\url{http://dx.doi.org/10.1007/BF00773456}
	
	\bibitem[{Gomes(2016)}]{AG-2016}
	Gomes, A.~J., 2016. Geometria extr\'inseca de campos de vetores em
	$\mathbb{R}^3$. http://repositorio.bc.ufg.br/tede/handle/tede/8636.
	
	\bibitem[{Gutierrez and Sotomayor(1991)}]{JSCG-1991}
	Gutierrez, C., Sotomayor, J., 1991. Lines of Curvature and Umbilic Points on
	Surfaces. $18^{th}$ Brasilian Math. Colloquium. IMPA.
	
	\bibitem[{Hirsch(1973)}]{HM-1973}
	Hirsch, M.~W., 1973. Stability of compact leaves of foliations. Dynamical
	systems, {\rm (}Proc. Sympos., Univ. Bahia, Salvador, 1971{\rm)}, 135--153.
	
	\bibitem[{Hirsch(1976)}]{HM-1976}
	Hirsch, M.~W., 1976. Differential topology. Springer-Verlag, New
	York-Heidelberg, graduate Texts in Mathematics, No. 33.
	
	\bibitem[{Hirsch et~al.(1977)Hirsch, Pugh, and Shub}]{HPS-1977}
	Hirsch, M.~W., Pugh, C.~C., Shub, M., 1977. Invariant manifolds. Lecture Notes
	in Mathematics, Vol. 583. Springer-Verlag, Berlin-New York.
	
	\bibitem[{Krouglov(2008)}]{Kr-2008}
	Krouglov, V., 2008. The curvature of contact structures on 3-manifolds. Algebr.
	Geom. Topol. 8~(3), 1567--1579.
	\newline\urlprefix\url{https://doi.org/10.2140/agt.2008.8.1567}
	
	\bibitem[{Lopes et~al.(2015)Lopes, Sotomayor, and Garcia}]{drs-2015}
	Lopes, D., Sotomayor, J., Garcia, R., 2015. Partially umbilic singularities of
	hypersurfaces of {$\mathbb{R}^4$}. Bull. Sci. Math. 139~(4), 431--472.
	\newline\urlprefix\url{https://doi.org/10.1016/j.bulsci.2014.10.005}
	
	\bibitem[{Monge(1796)}]{Mo-1796}
	Monge, G., 1796. Sur les lignes de courbure de la surface de l'ellipsoide.
	Journ. Ecole Polytech. II cah., 145--165.
	
	\bibitem[{Palis and de~Melo(1982)}]{PM-1982}
	Palis, Jr., J., de~Melo, W., 1982. Geometric theory of dynamical systems.
	Springer-Verlag, New York-Berlin, an introduction, Translated from the
	Portuguese by A. K. Manning.
	
	\bibitem[{Peixoto(1967)}]{MP-1967}
	Peixoto, M.~M., 1967. On an approximation theorem of {K}upka and {S}male. J.
	Differential Equations 3, 214--227.
	\newline\urlprefix\url{https://doi.org/10.1016/0022-0396(67)90026-5}
	
	\bibitem[{Reinhart(1977)}]{Re-1977}
	Reinhart, B.~L., 1977. The second fundamental form of a plane field. J.
	Differential Geom. 12~(4), 619--627.
	\newline\urlprefix\url{http://projecteuclid.org/euclid.jdg/1214434230}
	
	\bibitem[{Rifford(2014)}]{Rif-2014}
	Rifford, L., 2014. Sub-{R}iemannian geometry and optimal transport.
	SpringerBriefs in Mathematics. Springer, Cham.
	\newline\urlprefix\url{http://dx.doi.org/10.1007/978-3-319-04804-8}
	
	\bibitem[{Rifford and Ruggiero(2012)}]{RR-2012}
	Rifford, L., Ruggiero, R.~O., 2012. Generic properties of closed orbits of
	{H}amiltonian flows from {M}a\~n\'e's viewpoint. Int. Math. Res. Not.
	IMRN~(22), 5246--5265.
	
	\bibitem[{Rogers(1911/1912)}]{ro-1912}
	Rogers, G., 1911/1912. Some differential properties of the orthogonal
	trajectories of a congruence of curves. Proceedings of the Royal Irish
	Academy. Section A: Mathematical and Physical Sciences 29~(6), 92--117.
	
	\bibitem[{Rovenski and Walczak(2011)}]{RoWa-2011}
	Rovenski, V., Walczak, P., 2011. Topics in extrinsic geometry of
	codimension-one foliations. SpringerBriefs in Mathematics. Springer, New
	York, with a foreword by Izu Vaisman.
	\newline\urlprefix\url{https://doi.org/10.1007/978-1-4419-9908-5}
	
	\bibitem[{S\'{a}nchez-Bringas and Ram\'{\i}rez-Galarza(1995)}]{SB-1995}
	S\'{a}nchez-Bringas, F., Ram\'{\i}rez-Galarza, A.~I., 1995. Lines of curvature
	near umbilical points on surfaces immersed in {${\bf R}^4$}. Ann. Global
	Anal. Geom. 13~(2), 129--140.
	\newline\urlprefix\url{https://doi.org/10.1007/BF01120328}
	
	\bibitem[{Sotomayor and Gutierrez(1982)}]{CGJS-1982}
	Sotomayor, J., Gutierrez, C., 1982. Structurally stable configurations of lines
	of principal curvature. In: Bifurcation, ergodic theory and applications
	({D}ijon, 1981). Vol.~98 of Ast\'erisque. Soc. Math. France, Paris, pp.
	195--215.
	
	\bibitem[{Spivak(1999{\natexlab{a}})}]{MS01}
	Spivak, M., 1999{\natexlab{a}}. A comprehensive introduction to differential
	geometry V.1, 3rd Edition. Publish or Perish, INC.
	
	\bibitem[{Spivak(1999{\natexlab{b}})}]{MS03}
	Spivak, M., 1999{\natexlab{b}}. A comprehensive introduction to differential
	geometry V.3, 3rd Edition. Publish or Perish, INC.
	
	\bibitem[{Voss(1880)}]{Voss}
	Voss, A., 1880. Geometrische {I}nterpretation der {D}ifferentialgleichung
	{$Pdx+Qdy+Rdz=0$}. Math. Ann. 16~(4), 556--559.
	\newline\urlprefix\url{http://dx.doi.org/10.1007/BF01446224}
	
\end{thebibliography}
\end{document}